\documentclass[12pt]{article}

\usepackage{amsmath}
\usepackage{amsthm}
\usepackage{amsfonts}
\usepackage{stmaryrd}
\usepackage{amssymb}
\usepackage[all]{xy}
\usepackage{comment}

\newcommand\R{\mathbb{R}}

\newcommand\Z{\mathbb{Z}}

\newcommand\del{\partial}

\renewcommand{\cal}[1]{\mathcal{#1}}
\newcommand{\an}[1]{\langle{#1}\rangle}

\newcommand{\xra}{\xrightarrow}

\newcommand{\cvec}[2]{\left( \begin{array}{c} #1 \\ #2 \end{array} \right)}
\newcommand{\twotwo}[4]{\left( \begin{array}{cc} #1 & #2 \\ #3 & #4 \end{array} \right)}

\theoremstyle{plain}
\newtheorem{theorem}{Theorem}[section]
\newtheorem{lemma}[theorem]{Lemma}
\newtheorem{corollary}[theorem]{Corollary}
\newtheorem{proposition}[theorem]{Proposition}
\newtheorem{conjecture}[theorem]{Conjecture}
\newtheorem{question}[theorem]{Question}

\theoremstyle{definition}
\newtheorem{definition}[theorem]{Definition}
\newtheorem{example}[theorem]{Example}

\theoremstyle{remark}
\newtheorem{remark}[theorem]{Remark}

\usepackage{color}

\advance \topmargin by -\headheight
\advance \topmargin by -\headsep
     
\evensidemargin \oddsidemargin
\title{Intermediate curvatures and highly connected manifolds}
\author{Diarmuid Crowley and David J. Wraith}
\date{\today}

\begin{document}
\maketitle













\abstract{We show that after forming a connected sum with a homotopy sphere, 
all $(2j{-}1)$-connected $2j$-parallelisable manifolds in dimension $4j{+}1$,
$j \ge 2$, 
can be equipped with Riemannian metrics of 2-positive Ricci curvature. 
The condition of 2-positive Ricci curvature is defined to mean that the sum of the two smallest eigenvalues of the Ricci tensor is positive at every point. This result is a counterpart to a previous result of the authors concerning the existence of positive Ricci curvature on highly connected manifolds in dimensions $4j{-}1$ for $j \ge 2$, and in dimensions $4j{+}1$ for $j \ge 1$ with torsion-free cohomology. A key technical innovation involves performing surgery on links of spheres within 2-positive Ricci curvature.}

\section{Introduction} 
\label{sec:intro}
\leftskip 0cm \rightskip 0cm 
We begin by recalling that the term `highly connected manifold' means a manifold $M$ of dimension $2s$ or $2s{+}1$ which is $(s{-}1)$-connected, that is, for which all homotopy groups up to and including $\pi_{s-1}(M)$ vanish. We will assume without further comment that all manifolds in this paper are closed, connected and oriented unless stated otherwise. In [CW] the authors 
established the result below for positive Ricci curvature on highly connected manifolds in dimensions $4j{-}1$, and highly connected manifolds in dimensions $4j{+}1$ with torsion-free cohomology. We recall that a manifold is said to be $l$-parallelisable if its tangent bundle restricted to some $l$-skeleton is trivial.
\begin{theorem} [\cite{CW}; Theorems A and D] \label{thm:JDG}
Consider a $(2j{-}2)$-connected manifold $M^{4j-1}$ for $j \ge 2.$ If $j \equiv 1 \hbox{ mod }4$ assume further that $M$ is $(2j{-}1)$-parallelisable. Then there is a homotopy sphere $\Sigma^{4j-1}$ such that 
$M \sharp \Sigma$ admits a metric of positive Ricci curvature.
If $N^{4j+1},$ $j \ge 1,$ is a $(2j{-}1)$-connected $2j$-parallelisable manifold with torsion-free (integer) cohomology, then there is a homotopy sphere $\tilde{\Sigma}^{4j+1}$ such that $N \sharp \tilde{\Sigma}$ admits a Ricci positive metric.
\end{theorem}

The question of which manifolds admit metrics of positive Ricci curvature constitutes a major problem in Riemannian geometry. One might hope that some form of topological `simplicity' might facilitate the existence of Ricci positive metrics, and the above theorem illustrates that being highly connected in dimensions $4j{-}1$ is an appropriate interpretation of this term. Note that we cannot in general hope to remove the ambiguity created by having to form connected sums with homotopy spheres. It has long been known (\cite{Hi}) that there exist homotopy spheres which do not even satisfy the weaker positive scalar curvature condition, so there are highly-connected manifolds which definitely do not admit positive Ricci curvature.

It is natural to look for a counterpart to Theorem \ref{thm:JDG} involving highly-connected manifolds $M$ in dimensions $4j{+}1$ which goes beyond the highly-restrictive torsion-free cohomology requirement in that theorem.
Unfortunately the techniques used to establish Theorem 1.1 
do not apply when $H_{2j}(M^{4j+1})$ contains torsion.
(Here and throughout integer coefficients for (co)homology are omitted.)
The existence of positive Ricci curvature metrics on the manifolds of dimension $4j{-}1$ in Theorem 1.1 follows from the fact that one can express these manifolds as the boundaries of plumbed manifolds, where all the plumbings involve $D^{2j}$-bundles over $S^{2j}.$ The relevance of the plumbing construction to issues of positive Ricci curvature stems from the fact that on the boundary, the effect of plumbing two disc bundles over spheres is precisely a surgery on a fibre sphere of one of the bundles. Under limited circumstances (see \cite{Wr2} for example) one can extend a Ricci positive metric over a surgery, and crucially these circumstances are fulfilled in the case of surgeries arising from such plumbings. Now one can also use plumbings of this type to construct manifolds in dimensions $4j{+}1$, but it turns out that these manifolds belong to a short list of well-known objects, all having torsion-free integer cohomology: spheres, products or spheres and connected sums of products of spheres. These manifolds are all known to admit Ricci positive metrics (see for example \cite{SY}, \cite{Wr3}, \cite{BG3}).

Even though plumbing does not help in dimensions $4j{+}1$ when torsion is present in $H_{2j}(M)$, one can still hope to make progress using more general kinds of surgery. To this end we establish the following purely topological result:

\begin{theorem} \label{thm:surgery} 
Let $M^{4j+1}$ be a $(2j{-}1)$-connected $2j$-parallelisable manifold.
Then there is a homotopy sphere $\Sigma^{4j+1}$, a natural number $r \geq 0,$
and if $r > 0$ integers $m_1,...,m_r$ with each $m_i \geq 2$,
such that $M$ 
is diffeomorphic to a connected sum 
$$M_0 \sharp \Sigma \text{~($r = 0$)} \quad \text{or} \quad  
M_0 \sharp M_1 \sharp \cdots \sharp M_r \sharp \Sigma
\text{~($r > 0$)}.$$ 
Here $M_0$ is a connected sum of the total spaces of linear $S^{2j}$-bundles over $S^{2j+1}$ or possibly the empty connected sum $M_0 = S^{4j+1}$, and each $M_i$ for $i>0$ is constructed by performing surgery on a link of two embedded 
$2j$-spheres in $S^{4j+1}$ for which the linking number is equal to 
$m_i.$
\end{theorem}

We remark that the linking number referred to in the above theorem is the oriented intersection number between one of the $2j$-spheres and a disc $D^{2j+1}$ bounding the other sphere (see \cite[p.\,288]{ST} for more details). 

We also remark that the manifolds appearing in Theorem \ref{thm:surgery} we classified up to
connected sum with homotopy spheres by Wall \cite[Theorem 7]{Wa2}; the novelty of this theorem
is the explicit description of the manifolds $M_i$.
When $j = 1$ or $3$ the manifolds $M_i$ for $i>0$ have a long history in differential topology, featuring in Smale's early classification theorems for handlebodies and their boundaries \cite[\S 6]{Sm}.
(Note that Smale uses the notation $M_{m_i}$ for $M_i$.)

Turning our attention to geometry, the dimension $4j{+}1$ analogue of the dimension $4j{-}1$ statement in Theorem \ref{thm:JDG} is unfortunately still out of reach. This is a consequence of Theorem \ref{thm:surgery}: even though this is best topological description of highly-connected manifolds in dimension $4j{+}1$ at our disposal, we are nonetheless presented with a decomposition involving connected sums. The connected sum construction presents a problem for positive Ricci curvature. In the presence of non-trivial fundamental groups, it follows from the classical Bonnet-Myers Theorem that connected sums of Ricci positive manifolds do not in general admit positive Ricci curvature. This obstruction disappears in the simply-connected case, however little is known in this situation about the conditions under which Ricci positivity can be extended across connected sums. There are many examples of Ricci positive connected sums in the literature, but in all cases there is an alternative description of the manifolds which {\it avoids} the use of connected sums.

On the other hand, if we are prepared to weaken the curvature condition slightly, it turns out that we can indeed prove a counterpart theorem. The curvature condition in question is that the sum of the two smallest eigenvalues of the Ricci tensor is everywhere positive. (Of course for positive Ricci curvature we need all the eigenvalues positive, so this new condition allows at most one eigenvalue to be negative at any point, but with smaller magnitude than the smallest of the positive eigenvalues.) Following Wolfson (\cite{Wo}) we make the following definition:
\begin{definition} We say that an $n$-dimensional Riemannian manifold $(N,g)$ has 
$k$-positive Ricci curvature if the sum of the $k$ smallest eigenvalues of the Ricci tensor is positive at all points. We will write this as $Ric_k>0.$
\end{definition}
Thus we see that $n$-positive Ricci curvature is just positive scalar curvature, and 1-positive Ricci curvature is the same as positive Ricci curvature. The condition we will be interested in is 2-positive Ricci curvature. 

We note that a recent result of Wolfson, \cite{Wo2}, establishes that closed manifolds which admit metrics of 2-positive Ricci curvature have virtually free fundamental groups. This contrasts with the case of positive Ricci curvature, where by the Bonnet-Myers Theorem fundamental groups must be finite. 

The main aim of this paper is to establish:
\begin{theorem}\label{thm:main} 
Let $M^{4j+1}$ be a $(2j{-}1)$-connected $2j$-parallelisable manifold, $j \ge 2$.
Then there is a homotopy sphere $\Sigma^{4j+1}$ such that $M \sharp \Sigma$ admits a metric of 2-positive Ricci curvature.
\end{theorem}

Given the topological description presented in Theorem \ref{thm:surgery}, we remark that it is well known that linear $S^{2j}$-bundles over $S^{2j+1}$ admit metrics of positive Ricci curvature (see for example \cite[Corollary 3.6]{Na}). Moreover connected sums of such bundles also admit Ricci positive metrics by \cite[Theorem 2.4]{CW} in conjunction with \cite[Theorem 2.2]{Wr1}. Thus from a metric perspective, it is clear from Theorem \ref{thm:surgery} that Theorem \ref{thm:main} will follow if we can perform surgeries on links of $2j$-spheres in $S^{4j+1}$ within 2-positive Ricci curvature, and if we can perform connected sums preserving 2-positive Ricci curvature. This second point is contained in the following surgery result due to Wolfson (\cite{Wo}, but see also \cite{Ho}). To the best of our knowledge, this is the only result to date in the literature concerning the existence of 2-positive Ricci curvature metrics, (but see \cite{WW} and \cite{Kor} which consider spaces of metrics).

\begin{theorem}[\cite{Wo},\cite{Ho}]  \label{thm:wolfson} 
Let $M^n$ be a closed Riemannian manifold with $k$-positive Ricci curvature, $2 \le k \le n$. Then any manifold obtained from $M$ by performing surgeries in codimension $q$ with $q \ge \max\{n+2-k,3\}$ also admits a metric of $k$-positive Ricci curvature. In particular if $M_1$ and $M_2$ are manifolds of dimension $n \ge 3$ which admit metrics of 2-positive Ricci curvature, then the connected sum $M_1 \sharp M_2$ also admits a metric of 2-positive Ricci curvature.
\end{theorem}

It is therefore clear that our main metric focus in this paper will be to establish the existence of 2-positive Ricci curvature metrics on the manifolds $M_i$, $i>0$, appearing in Theorem \ref{thm:surgery}. Performing the necessary metric surgeries presents considerable technical difficulties, and finding an approach to overcoming these problems is perhaps the most novel aspect of this paper. The key idea is to create a metric on the $(4j+1)$-sphere which is adapted to the topology of the embedded $2j$-spheres on which we will perform surgery.
 
A major difference between surgery preserving positive scalar curvature and surgery preserving $Ric_2>0$ (or positive Ricci curvature) is that the former surgery is local in nature, whereas the latter is not. The known Ricci positive surgery results require large normal bundles compared to the size of the surgery sphere. In contrast, the obvious picture of linked spheres involves small normal bundles compared to the size of the links. Indeed one of the principal problems in extending the scope of Ricci positive surgery theory is understanding when and how one can handle surgery on linked spheres. Very little is known about this at present. However, the ideas used in proving Theorem \ref{thm:main} for the weaker $Ric_2>0$ condition might offer a blueprint for making progress in positive Ricci curvature.  

The techniques we use to establish Theorem \ref{thm:main} unfortunately fail in dimension five. It seems likely, however, that with some adaptation of the constructions we make, it might be possible to extend our results to include simply-connected five-manifolds. We therefore propose 

\begin{conjecture}\label{con:dim5} 
Every simply-connected 5-manifold admits a metric of 2-positive Ricci curvature.
\end{conjecture}

Conjecture \ref{con:dim5} is of interest as there is a long-standing open question as to whether all simply-connected 5-manifolds admit metrics of positive Ricci curvature. There is a classification simply-connected 5-manifolds due to Barden (\cite{Ba}, but see also \cite{Sm} for the spin case). These manifolds fall into two infinite families depending on whether they are spin or non-spin. In the spin case, most of the manifolds were shown to admit Ricci positive metrics by the work of Boyer and Galicki on Sasakian geometry (\cite{BG1}, \cite{BG2}), though this existence result still omits inifinitely many objects in the spin class. In contrast, until very recently only two of the non-spin 5-manifolds have been known to admit Ricci positive metrics, though this set has now been expanded in \cite{CGG}.

We also note that in the analogous - though much stronger - context of the curvature operator of a Riemannian manifold, the 2-positive condition has received much attention, culminating in the work of B\"ohm and Wilking \cite{BW}.
\bigskip

By summing the $k$ smallest eigenvalues of the Ricci tensor we obtain in a simple way a natural family of curvatures which interpolate between scalar curvature and Ricci curvature. This raises an obvious question: given a phenomenon which holds for positive scalar curvature say, but which might not hold for positive Ricci curvature, for which of these intermediate curvatures does the phenomenon hold? The analysis of this question clearly has the potential to offer a deeper insight into properties observed in the more classical worlds of scalar and Ricci curvatures. In the simply-connected case, there is no known difference between the class of closed positive scalar curvature manifolds and closed Ricci positive manifolds in dimensions at least 5, and therefore no known difference between these classes and the class of $k$-positive Ricci curvature manifolds for any $k$. In fact it was shown in (\cite{Wo}; Theorem 2.3) that the class of closed simply-connected $n$-manifolds ($n \ge 5$) which admit positive scalar curvature is precisely the same as the corresponding class of $(n{-}1)$-positive Ricci curvature manifolds. On the other hand there are differences between positive scalar curvature and $k$-positive Ricci curvature for $1 \le k <n$ in the non-simply-connected case. To see this consider the Riemannian product manifold $S^m(R) \times (M^k,g)$ where the first factor is a round sphere of radius $R$ and dimension $m \ge 2$, and the second factor is a closed hyperbolic spin manifold of dimension $k \ge 2$. Then for $R$ sufficiently small this product manifold has $(k+1)$-positive Ricci curvature, as any collection of $k+1$ linearly independent tangent directions must include a direction with a component in $TS^m.$ However this manifold does not support a metric of $k$-positive Ricci curvature, since if it did, this would mean that $M^k$ admits a positive scalar curvature metric, as this is not possible by \cite{GL2}.

In the context of the above discussion, it seems natural to pose the following
\begin{question}\label{q}
To what extent do the classes of simply-connected closed manifolds admitting positive scalar curvature metrics and 2-positive Ricci curvature metrics differ? In particular, do these two classes agree?
\end{question}

Note that Theorem \ref{thm:3c4p9m} answers the latter question affirmatively for $3$-connected $4$-parallelisable $9$-manifolds.
\bigskip

In the same way that there is a natural gradation of curvatures between scalar and Ricci, there are similarly natural notions of curvature intermediate between the Ricci and sectional curvatures. If we fix a tangent vector $u$ and extend it to an orthogonal $(k+1)$-frame $\{u,v_1,...,v_k\}$, we can form the sum of sectional curvatures $\sum_{i=1}^k K(u,v_i).$ A manifold is said to have $k$-positive curvature if these sums are all positive. This condition has been studied much more than $k$-positive Ricci curvature: see for example \cite{Ha}, \cite{Sh1}, \cite{Sh2}, \cite{Wu}, \cite{Wi}, \cite{GX}, \cite{GW1}, \cite{GW2}, \cite{GW3}, \cite{Mo1}, \cite{Mo2}.

We note, however, that the terminology for these intermediate curvature notions is not fixed in the literature. For example in \cite{GW1}, \cite{GW2}, \cite{GW3} what we have called $k$-positive curvature above is referred to as the $k^{th}$-intermediate Ricci curvature, and denoted $Ric_k$. In \cite{Wi} it is called the $k^{th}$-Ricci curvature. There is an obvious and significant potential for confusion here.

We should also mention that there are other notions of curvature which are in some sense intermediate between classical curvatures. These include Labbi's $p$-curvature, see for example \cite{La}, \cite{BL1}, and the $\sigma_k$ and $\Gamma_k$ curvatures defined in terms of the eigenvalues of the Schouten tensor, see for example \cite{BL2}.
\bigskip

This paper is laid out as follows. In \S2 we study the topology of highly connected manifolds in dimensions $4j{+}1$, leading to a proof of Theorem \ref{thm:surgery}. In \S3 we show how to embed pairs of $2j$-spheres into $S^{4j+1}$ in a geometrically nice way, so as to realize any given linking number. In \S4 we show how to construct a $Ric_2>0$ metric on $S^{4j+1}$ which facilitates our final step, which is to perform 2-positive Ricci curvature surgery on the link. This is carried out in \S5, concluding with the proof of Theorem \ref{thm:main}. There is also an Appendix which explains some of the more intricate curvature formulas appearing in the main body of the paper.
\smallskip

\noindent {\bf Acknowlegement:} The second author would like to thank the University of Aberdeen, where this work was started, for their hospitality.
%

\section{$(2j{-}1)$-connected $(4j{+}1)$-manifolds} \label{sec:top}

\subsection{The general case} \label{subsec:top-gen}
To state our main topological result, Theorem \ref{thm:top_main} below,
we first recall some notation from \cite{CW}.
Let $\pi_m^s = \mathrm{lim}_{i \to \infty} \pi_{m+i}(S^i)$ be 
the stable $m$-stem and $J_m \colon \pi_m(SO) \to \pi_m^s$
the $J$-homomorphim.  By \cite{KM} there is an exact sequence
\[ 0 \to bP_{4j+2} \to \Theta_{4j+1} \xra{~\Phi~} \mathrm{Coker}(J_{4j+1}) \to 0, \]
where $\Theta_{4j+1}$ is the group of oriented diffeomorphism classes
of homotopy $(4j{+}1)$-spheres and $bP_{4j+2} \subset \Theta_{4j+1}$
is the subgroup of those homotopy spheres which bound parallelisable manifolds.

Recall that $M$ is a $(2j{-}1)$-connected $(4j{+}1)$-manifold.
For the majority of this section we also assume that the tangent 
bundle of $M$ is trivial over any $2j$-skeleton of $M$;
i.e.~that $M$ is $2j$-parallelisable.
In \cite[\S 5.1]{CW} we recalled 
the bordism group $\Omega_{4j+1}^{O\an{2j}}$ of $2j$-parallelisable
manifolds and showed that an orientation on $M$ determines
a bordism class $[M, \bar \nu] \in \Omega_{4j+1}^{O\an{2j}}$, where 
$\bar{\nu}:M \to BO\an{2j{+}1}$ is a lift of the map $M \to BO$ which classifies the stable normal bundle of $M$,
and $BO\an{2j{+}1} \to BO$ is the $2j$-connected covering of $BO$. 
By \cite[Theorem 7.1]{CW} there is a homotopy sphere $\Sigma$ such that
$[M, \bar \nu] = [\Sigma] \in \Omega_{4j+1}^{O\an{2j}}$.  It follows that
the natural map $\eta \colon \mathrm{Coker}(J_{4j+1}) \to \Omega_{4j+1}^{O\an{2j}}$
is onto and we define
\[ \Sigma(M) \subset \Theta_{4j+1} \]
to be the set of  homotopy spheres $\Sigma$ such that 
$\eta(\Phi(\Sigma)) = - [M, \bar \nu] \in \Omega_{4j+1}^{O\an{2j}}$. 


\begin{theorem}  \label{thm:top_main}
Let $M^{4j+1}$ be a $(2j{-}1)$-connected $2j$-parallelisable manifold.
Then there is a natural number $r \geq 0,$ and if $r > 0$ 
a collection of integers $m_1,...,m_r$
with each $m_i \geq 2$ such that 
for every homotopy sphere $\Sigma \in \Sigma(M)$,
$M \# \Sigma$ is diffeomorphic to a manifold $M_0$ (if $r=0$) or to a connected sum 
%
$$M_0 \sharp M_1 \sharp \cdots \sharp M_r,$$
%
where $M_0$ is a connected sum of the total spaces of linear 
$S^{2j}$-bundles over $S^{2j+1}$ or possibly the empty connected sum $M_0 = S^{4j+1}$,
and each $M_i$ for $i>0$ is constructed by performing surgery on a two component link of embedded $2j$-spheres in $S^{4j+1}$ for which the linking number is equal to $m_i.$
\end{theorem}

Regarding the linking numbers $m_i$, we can restrict to the case $m_i \geq 2$ for
the following reasons:
if $m_i = 1$, then the outcome
of surgery on the corresponding link is a homotopy sphere, if $m_i = 0$, then the outcome
of surgery is diffeomorphic to the connected sum of two $S^{2j}$-bundles over $S^{2j+1}$
and if $m_i < 0$, then the outcome of surgery is diffeomorphic to the outcome of surgery
on $-m_i$.  As discussed in the introduction, the manifolds described in 
Theorem \ref{thm:top_main} were classified up to connected sum with homotopy
spheres in \cite[Theorem 7]{Wa2}.

To give the proof of Theorem \ref{thm:top_main} we must first recall
foundational results of Smale and Wall on certain handlebodies.  
We have already summarised these results in [CW \S4] and we use the notation
and setting we developed there.  We define a {\it handlebody} $W$ of dimension $4j{+}2$ to 
be a smooth manifold obtained by attaching $(2j{+}1)$-handles to the 
$(4j{+}2)$-disc $D^{4j+2}$, so that
\begin{equation} \label{eq:hb}
W = D^{4j+2} \cup_{\phi} (\sqcup_{i=1}^s D^{2j+1} \times D^{2j+1}),
\end{equation}
where $\phi \colon \sqcup_{i=1}^s D^{2j+1} \times S^{2j} \to S^{4j+1}$ is a 
smooth embedding.
We define
$$  \mathcal{H}^{4j+2} := \{ W\,|\, W \cong  D^{4j+2} \cup_{\phi} (\sqcup_{i=1}^s
D^{2j+1} \times D^{2j+1}) \} $$
to be the set of 
diffeomorphism classes of 
handlebodies.  

We next review Wall's classification of $\mathcal{H}^{4j+2}$ \cite{Wa1}.
Recall from [CW \S3] the {\it quadratic form parameter} 
$$\pi_{2j}\{SO(2j{+}1)\} = \bigl( \pi_{2j}(SO(2j{+}1)), \Z, {\rm h}, {\rm p} \bigr).$$
An extended quadratic form with values in $\pi_{2j}\{SO(2j{+}1)\}$
is a triple $(H, \lambda, \mu)$ where $H$ is a finitely generated free abelian group,
$\lambda \colon H \times H \to \Z$ is a skew-symmetric bilinear form and 
$\mu \colon H \to \pi_{2j}(SO(2j{+}1))$ is a function such that
$$ \mu(x + y) = \mu(x) + \mu(y) + {\rm p}(\lambda(x, y)) 
\quad \hbox{and} \quad {\rm h}(\mu(x)) = \lambda(x, x). $$
There are obvious notions of isometry and orthogonal sum for extended quadratic forms.
Given a handlebody $W$ there is a well-defined function,
$$ \mu_W \colon H^{2j+1} (W,\partial W) \to \pi_{2j}(SO(2j{+}1)), \quad
x \mapsto \nu_{\hat x},$$
which is defined by taking the isomorphism class of the normal bundle $\nu_{\widehat x}$ 
of an embedding
$\widehat{x} \colon S^{2j+1} \to W$ which represents the Poincar\'{e} dual of $x$.
Moreover, Wall showed that if $\lambda_W \colon H^{2j+1} (W, \partial W) \times H^{2j+1} (W, \partial W) \to \Z$ is the intersection form of $W$, then the triple 
$$(H^{2j+1} (W, \partial W), \lambda_W, \mu_W)$$
defines an extended intersection form with values in $\pi_{2j}\{SO(2j{+}1)\}$ and proved the following 

\begin{theorem}[{\cite[p.\,168]{Wa1}}] \label{thm:wall_4j+2_classification}
For all $j \geq 1$,
the assignment of its extended intersection form to a handlebody defines a bijection,
\[  \cal{H}^{4j+2} \equiv \cal{F}^{4j+2}, \quad W \mapsto (H^{2j+1} (W, \del W), \lambda_W, \mu_W), \]
which maps the boundary connected sum of handlebodies to the orthogonal sum of forms;
\[ W_0 \natural W_1 \mapsto (H^{2j+1} (W_0, \del W_0), \lambda_{W_0}, \mu_{W_0}) \oplus 
(H^{2j+1} (W_1, \del W_1), \lambda_{W_1}, \mu_{W_1}).\]
Moreover, every isomorphism of extended intersection forms,
\[ A \colon (H^{2j+1} (W_1, \del W_1), \lambda_{W_1}, \mu_{W_1}) \cong 
(H^{2j+1} (W_0, \del W_0), \lambda_{W_0}, \mu_{W_0}), \]
is realised by a diffeomorphism $f_A \colon W_0 \cong W_1$.
\end{theorem}

\begin{example} \label{ex:1link}
If $W$ is obtained from $D^{4j+2}$ by attaching a single $(2j{+}1)$-handle,
so that $s = 1$ in \eqref{eq:hb}, then the embedding 
$\phi \colon D^{2j+1} \times S^{2j} \to S^{4j+1}$ is
a framed $2q$-sphere which has self-linking number $0$.
In this case the extended quadratic form of $W$ is given by
\[ \bigl( H^{2j+1} (W, \del W), \lambda_W, \mu_W \bigr) \cong
\left(\Z(x_1), 0, \gamma_1 \right), \]
where 
$\mu_W(x_1) = \gamma_1 \in \pi_{2j}(SO(2j{+}1))$.
\end{example}

\begin{example} \label{ex:2link}
If $W$ is obtained from $D^{4j+2}$ by attaching two $(2j{+}1)$-handles,
so that $s = 2$ in \eqref{eq:hb}, then the embedding 
$\phi \colon (D^{2j+1} \times S^{2j}_1) \sqcup (D^{2j+1} \times S^{2j}_2) \to S^{4j+1}$ is
a framed link with two components which has a linking number $m \in \Z$.
In this case the extended quadratic form of $W$ is given by
\[ \bigl( H^{2j+1} (W, \del W), \lambda_W, \mu_W \bigr) \cong
\left(\Z^2(x_1, x_2), \twotwo{0}{m}{-m}{0} , \cvec{\gamma_1}{\gamma_2} \right), \]
where 
$\mu_W(x_i) = \gamma_i \in \pi_{2j}(SO(2j{+}1))$.
\end{example}

\begin{proof}[Proof of Theorem \ref{thm:top_main}]
By \cite[Theorem 7.1]{CW} $M \# \Sigma$ is diffeomorphic to the boundary
of a handlebody $W$.  For brevity let
$(H, \lambda, \mu) : = (H^{2j+1} (W, \del W), \lambda_W, \mu_W)$
be the extended intersection form of $W$
and define {\em the radical of $\lambda$}, $H_{0} \subset H$, by
$$H_{0} := \{x \in H \, | \,\lambda(x, y) = 0~\forall y \in H \}.$$
There is an orthogonal decomposition
$(H, \lambda, \mu) = (H_{0}, 0, \mu_{0}) \oplus (H_{t}, \lambda_t, \mu_t)$,
where $H_t \subset H$ is a complementary summand to $H_{0}$ and
$\lambda_t := \lambda|_{H_t \times H_t}$ is non-degenerate.
Note that it is possible that $H_t = 0$, in which case $r = 0$.
By the classification of non-degenerate skew symmetric forms, see 
\cite[\S 14]{VF},
there is an isometry 
$$(H_t, \lambda_t) \cong \bigoplus_{i=1}^r (H_{i}, \lambda_{i}), $$
where for every $i$, there is a positive integer $m_i$ and an isometry
$$ (H_{i}, \lambda_{i}) \cong 
\left(\Z^2(x_{i1}, x_{i2}), \twotwo{0}{m_i}{-m_i}{0} \right).$$
It follows that there are isometries of extended quadratic forms
$$ (H_{t}, \lambda_{t}, \mu_t) \cong \bigoplus_{i=1}^r
\left(\Z^2(x_{i1}, x_{i2}), \twotwo{0}{m_i}{-m_i}{0}, \cvec{\gamma_{i1}}{\gamma_{i2}} \right)$$
and
$$ (H, \lambda, \mu) \cong 
\bigoplus_{l=1}^s \bigl(\Z(x_l), 0, \gamma_{l} \bigr) \oplus
\bigoplus_{i=1}^r
\left(\Z^2(x_{i1}, x_{i2}), \twotwo{0}{m_i}{-m_i}{0}, \cvec{\gamma_{i1}}{\gamma_{i2}} \right),$$
where $\{x_1,  \dots, x_s\}$ is a basis for $H_{0}$ and
we allow $s = 0$ or $r = 0$, in which case the corresponding summand is empty.
By Theorem \ref{thm:wall_4j+2_classification} and 
Examples \ref{ex:1link} and \ref{ex:2link} it follows that $W$ is either diffeomorphic to 
$D^{4j+1}$ or to a
boundary connected sum of the form
$$ W \cong  \natural_{l=1}^s W_l,
\quad \text{or} \quad
W \cong  \natural_{i=1}^r W_i,
\quad \text{or} \quad
W \cong \bigl( \natural_{l=1}^s W_l \bigr) \natural \bigl( \natural_{i=1}^r W_i \bigr),$$
where every $W_l$, $l = 1, \dots, s$, is the total space of a $D^{2j+1}$-bundle over $S^{2j+1}$
and every $W_i$, $i = 1, \dots, r$, is obtained by attaching two $(2j{+}1)$-handles
to $D^{4j+2}$.
Now the manifold $M_l := \del W_l$, $l = 1, \dots, s$, is the total space of an $S^{2j}$-bundle 
over $S^{2j+1}$ and 
the manifold $M_i := \del W_i$, $i = 1, \dots, r$, 
is constructed by performing surgery on a two component link of embedded $2j$-spheres in $S^{4j+1}$ for which the linking number is equal to $m_i$.
Since there is a diffeomorphism 
$M \# \Sigma \cong \del W = (\#_{l=1}^s M_l) \# (\#_{i=1}^r M_i)$, 
the theorem follows.
\end{proof}


\subsection{$3$-connected $4$-parallelisable   $9$-manifolds} \label{subsec:1c5m}
Let $M$ be a $3$-connected $9$-manifold.  
Up to equivalence, $M$ admits a unique spin structure,
and a spin $9$-manifold $(M, s)$ has an $\alpha$-invariant
$ \alpha(M, s) \in KO_9 = \Z/2$
defined by taking the index of the Dirac operator associated to $s$.
Hence $M$ has a well-defined $\alpha$-invariant
$$ \alpha(M) := \alpha(M, s)  \in KO_9.$$

\begin{theorem} \label{thm:3c4p9m}
Let $M$ be a $3$-connected $4$-parallelisable $9$-manifold.
Then the following are equivalent:
\begin{enumerate}
\item[(a)] $\alpha(M) = 0 \in KO_9$;
\item[(b)] $M$ admits a metric of positive scalar curvature;
\item[(c)] $M$ admits a metric of $2$-postive Ricci curvature.
\end{enumerate}
\end{theorem}

\begin{proof} 
We have that $(a)$ and $(b)$ are equivalent by \cite{St}. The implication from 
$(c)$ to $(b)$ is trivial. Finally we show that $(a)$ implies $(c)$.
The $\alpha$-invariant defines a homomorphism 
$\alpha \colon \mathrm{Coker}(J_9) \to KO_9$.
By \cite[Lemma 5.8]{BCS} there is a homotopy sphere $\Sigma$ for which
$\Phi(\Sigma)$ generates $\mathrm{ker}(\alpha)$ and which is the boundary of
a plumbing manifold $W$ which is obtained from $D^{10}$ by
adding a $4$-handle and a $6$-handle.
By \cite{Wr1}, $\Sigma$ admits a metric of positive Ricci curvature.
Since $2$-positive Ricci curvature can be preserved under connected sums,
the theorem now follows from Theorem \ref{thm:top_main}.
\end{proof}

%


\bigskip\medskip
\section{Linked spheres in $S^{2k+1}$}
\label{sec:links}
\bigskip\medskip

From now on we will work in dimension $2k+1$, assuming $k \ge 3$. Although the applications we have in mind concern dimensions $4k+1$, the topological and geometric constructions we make in this and subsequent sections generally work in dimension $2k+1$. We will therefore work in this more general setting. 

The aim of this section is to show how to construct a pair of linked $k$-spheres within $S^{2k+1}$ which realize a given linking number $m \in \Z$, and which facilitate the metric surgery constructions to be performed in the subsequent sections.

We will view the ambient sphere as a union
\begin{align*}
S^{2k+1}&=(S^k_a \times D^{k+1}_b)\cup(D^{k+1}_a \times S^k_b)\\
&=\Bigl(S^k_a \times \{0\}\Bigr)\cup \Bigl((0,T) \times S^k_a \times S^k_b\Bigr) \cup \Bigl(\{0\} \times S^k_b\Bigr) \\
\end{align*}
for some $T>0.$ The round metric on $S^{2k+1}$ can be neatly expressed as a double warped product using this last viewpoint: $$ds^2_{2k+1}=dt^2+\cos^2(t)ds^2_k+\sin^2(t)ds^2_k,$$ where $T=\pi/2,$ the warping factor $\cos(t)$ scales $S^k_a,$ and $\sin(t)$ scales $S^k_b.$ The metric on $S^{2k+1}$ that we will need in the next section will, in some sense, be a generalization of this double warped product idea, though the metric restricted to the $S_b^k$ will not always be round, and in certain directions the metric will be `twisted'.

Let us set $Z^k:=S^k_a \times \{0\}.$ We will call this the {\it core sphere}. This will be one part of the two-component link of spheres on which we will perform surgery. Choose a basepoint $\ast \in S^k_b$, and consider the disc $${\mathcal D}:=\Bigl(S^k_a \times \{0\}\Bigr) \cup \Bigl((0,T) \times S^k_a \times \{\ast\}\Bigr) \cup \Bigl(\{0\} \times \{\ast\}\Bigr).$$ As this clearly bounds the core sphere, we will call $\mathcal D$ the {\it core disc}.

Our main task in this section is to construct an embedded sphere $W^k \subset S^{2k+1}$ such that given any $m \in \Z$, we have linking number $Lk(W,Z)=m.$
\medskip

We begin by considering a smooth function $f:S^k_a \to S^k_b$ of degree $d.$ This gives rise to a graph $\Gamma(f) \subset S^k_a \times S^k_b,$ which of course is a smoothly embedded copy of $S^k$ within $S^k_a \times S^k_b.$ If we let $\pi:S^k_a \times S^k_b \to S^k_a$ denote the standard projection map, then $\Gamma(f)$ is the image of a section for this trivial bundle. 

\begin{lemma}
For any $x \in S^k_a,$ the submanifold $\{x\} \times S^k_b \subset S^k_a \times S^k_b$ is transverse to $\Gamma(f)$ at the point $(x,f(x)).$
\end{lemma}

\begin{proof}
Since $\Gamma(f)$ is a graph, it is evident that the derivative map $\pi_*$ maps the tangent space $T_{(x,f(x))}\Gamma(f)$ isomorphically onto $T_x S^k_a.$ Thus $T_{(x,f(x))}(\{x\} \times S^k_b)={\mathrm{ker}}\pi_*$ is a complementary subspace to $T_{(x,f(x))}\Gamma(f)$ within $T_{(x,f(x))}(S^k_a \times S^k_b),$ hence the result.
\end{proof}

For convenience we will temporarily introduce a background metric $g_0=ds^2_k+ds^2_k$ on $S^k_a \times S^k_b$. With respect to $g_0$ we can consider the normal bundle $\nu_{g_0}\Gamma(f)$ of $\Gamma(f)$ within $S^k_a \times S^k_b$. Of course the topological properties of the normal bundle are independent of the metric, and therefore in some of our topological considerations below it will make sense to suppress the metric from our notation.

Let $H_x$ denote the hemisphere within $\{x\} \times S^k_b$ (as determined by $g_0$) centered on the point $(x,f(x)),$ and
consider the submanifold $D_\perp(\Gamma(f)) \subset S^k_a \times S^k_b,$ defined by $$D_\perp(\Gamma(f)):=\{H_x\,|\,x \in S^k_a\}.$$ It is clear that there is a projection map $\hat\pi:D_\perp(\Gamma(f)) \to \Gamma(f)$ given by $(x,y) \mapsto (x,f(x)),$ making $D_\perp(\Gamma(f))$ into a disc bundle over $\Gamma(f)$. We immediately deduce
\begin{corollary}
\label{tubular}
A tubular neighbourhood of $\Gamma(f)$ has the structure of a disc-bundle over $\Gamma(f)$ with fibres $D_x^k \subset \{x\} \times S^k_b$ centered on $(x,f(x)).$
\end{corollary}

For each $x \in S^k_a,$ consider the tangent space to $H_x$ at the point $(x,f(x)).$ Collectively, these tangent spaces form a vector bundle over $\Gamma(f)$ which we will denote $E(\Gamma(f)).$ Observe that $D_\perp(\Gamma(f))$ can be obtained by exponentiating (with respect to $g_0$) a disc bundle $DE(\Gamma(f)).$ In particular, $D_\perp(\Gamma(f))$ and $DE(\Gamma(f))$ are equivalent fibre bundles.

It is easy to construct a metric on $S^k_a \times S^k_b$ for which the normal bundle at $\Gamma(f)$ is precisely $E(\Gamma(f))$, and thus we can take $E(\Gamma(f))$ as a model for $\nu(\Gamma(f))$. In the next section our aim will be to construct a particular metric of this kind, hence the choice of notation for the bundle $D_\perp(\Gamma(f))$.

\begin{lemma}
\label{Euler_number}
If $k$ is even, the normal bundle $\nu(\Gamma(f))$ has Euler number $2d$. Thus $D_\perp(\Gamma(f))$ is non-trivial whenever $k$ is even and $d \neq 0.$
\end{lemma}

\begin{proof}
We begin by observing that the graph $\Gamma(f)$ is clearly represented by the homology class $(1,d) \in H_k(S^k_a \times S^k_b) \cong \Z \oplus \Z.$

The desired Euler number is equal to the self-intersection number of $\Gamma(f)$ within $S^k_a \times S^k_b$. We claim that this self-intersection number is $2d$.

In order to compute the self-intersection, we switch from homology to cohomology and consider cup products. The Poincar\'e dual of the homology class representing $\Gamma(f)$ is $(1,d) \in H^k(S^k_a \times S^k_b),$ and we need to compute the cup product $(1,d) \cup (1,d) \in H^{2k}(S^k_a \times S^k_b) \cong \Z.$

By the K\"unneth theorem, as rings we have $H^*(S^k_a \times S^k_b) \cong H^*(S^k_a) \otimes H^*(S^k_b),$ and therefore 
\begin{align*}
H^{2k}(S^k_a \times S^k_b) &\cong H^k(S^k_a) \otimes H^k(S^k_b) \cong \Z \otimes \Z \cong \Z \,\, \hbox{ and}\\
H^k(S^k_a \times S^k_b) &\cong \Bigl(H^0(S^k_a) \otimes H^k(S^k_b)\Bigr) \oplus \Bigl(H^k(S^k_a) \otimes H^0(S^k_b)\Bigr). \\
\end{align*}

Let $\alpha,\beta$ denote generators of $H^k(S^k_a)$ respectively $H^k(S^k_b).$ Via the last K\"unneth isomorphism we can interpret $(1,d) \in H^k(S^k_a \times S^k_b)$ as $(1 \otimes d\beta)\oplus (\alpha \otimes 1).$ Therefore $(1,d)\cup (1,d)$ can be viewed as 
\begin{align*}
\Bigl(&(1 \otimes d\beta)\oplus (\alpha \otimes 1)\Bigr) \cup \Bigl((1 \otimes d\beta)\oplus (\alpha \otimes 1)\Bigr)\\
&=(1 \otimes (d\beta)^2) + (\alpha^2 \otimes 1) + d(1\otimes \beta)\cup(\alpha\otimes 1)+d(\alpha\otimes 1)\cup (1\otimes \beta)\\
&=2d(\alpha \otimes \beta).\\
\end{align*}
The last line follows as the cup product is symmetric as a consequence of $k$ being even, together with the fact that $\alpha^2=\beta^2=0$. Now $H^{2k}(S^k_a \times S^k_b)$ is generated by the element $\alpha \otimes \beta$, so denoting the fundamental homology class of $S^k_a \times S^k_b$ by $[S^k_a \times S^k_b]$ we conclude that $$\Bigl((1,d) \cup (1,d)\Bigr) \cap [S^k_a \times S^k_b]=2d,$$ which establishes the claim.
\end{proof}

For each point $(x,y) \in S^k_a \times S^k_b,$ the background metric $g_0$ allows us to unambiguously identify the `$S^k_b$-antipodal point' $(x,-y).$ Thus we obtain an `antipodal section' to $\Gamma(f)$ which we will denote $\hat{\Gamma}(f).$ We obtain a disc bundle $D_\perp(\hat{\Gamma}(f))$ over this antipodal section exactly as we did for $\Gamma(f).$ Clearly $D_\perp({\Gamma}(f))$ and $D_\perp(\hat{\Gamma}(f))$ are equivalent fibre bundles, which by Lemma \ref{Euler_number} are non-trivial whenever $k$ is even and $d \neq 0.$ We immediately deduce

\begin{corollary}
\label{decomp}
$$S^k_a \times S^k_b=D_\perp(\Gamma(f))\cup_{id}D_\perp(\hat\Gamma(f)).$$
\end{corollary}

Notice that we can `combine' the two disc bundles above to obtain a single bundle given by $\bar\pi:S^k_a \times S^k_b \to \Gamma(f)$, where $\bar\pi(x,y)=(x,f(x)).$ If we let $\sigma:\Gamma(f) \to S^k_a$ denote the diffeomorphism $\sigma(x,f(x))=x,$ then the composition $\sigma \circ \bar\pi$ clearly agrees with the standard projection map $\pi:S^k_a \times S^k_b \to S^k_a.$ Thus $\bar\pi$ defines a trivial bundle structure, despite being decomposable into a union of typically non-trivial bundles over the same base space.

\begin{definition} 
Let $W$ denote the $k$-sphere $$W:=\{T/2\} \times \Gamma(f) \subset (0,T) \times S^k_a \times S^k_b \subset S^{2k+1}.$$
\end{definition}

\begin{proposition}
\label{linking_number}
${\mathrm Lk}(W,Z)=d.$
\end{proposition}

\begin{proof} 
We begin by noting the following result (\cite[p.\,72]{Ko}): if $h:(M^n,\partial M) \to (N^n,\partial N)$ has degree $\delta,$ then $\delta=I(\Gamma(h),M_p),$ where $\Gamma(h)$ is the graph of $h$ in $M \times N,$ $M_p=M \times \{p\} \subset M \times N$ is any `$M$-slice' through the product, and $I$ denotes the oriented intersection number. In our situation this means that $d$ is the intersection number of $\Gamma(f) \subset S^k_a \times S^k_b$ with any $\{x\} \times S^k_b.$ Now observe that this intersection number is trivially also the intersection number of $W$ with the core disc, i.e. $d=I(W,{\mathcal D}).$ Finally, as ${\mathcal D}$ bounds the core sphere $Z$, by the definition of linking numbers (see for example \cite[p.\,288]{ST}) we have $d=I(W,{\mathcal D})={\mathrm Lk}(W,Z)$ as claimed.
\end{proof}

\begin{lemma}
\label{s-trivial}
The normal bundle $\nu(\Gamma(f))$ is stably-trivial.
\end{lemma}

\begin{proof}
As the topological properties of the any normal bundle are independent of the Riemannian metric, for simplicitly let us equip $S^{2k+1}$ with the round metric, which we can view as a double warped product as outlined at the start of this section. We then have a smooth one-parameter family of embeddings of $S^k$ into $S^{2k+1}$ given by $\{t\} \times \Gamma(f)$ for $t \in [0,\pi/4],$ which starts with $Z$ and ends with $W$. Notice that when $t=0$ the sphere $S^k_b$ collapses to a point, and hence the graph $\Gamma(f)$ should be interpreted as $Z=S^k_a \times \{0\}.$ Thus within $S^{2k+1},$ $Z$ and $W$ must have isomorphic normal bundles. Clearly the normal bundle of $Z$ is trivial, and thus it follows that the normal bundle of $W$ is also trivial. Finally, we observe that $\nu(\Gamma(f))\oplus \R$ agrees with the normal bundle of $W$ within $S^{2k+1}.$
\end{proof}


\bigskip\medskip
\section{Metrics on $S^{2k+1}$}
\label{sec:metrics}
\bigskip\medskip


The aim of this section is to define a metric on $S^{2k+1}$ with $Ric_2>0$ which will allow surgery on $Z$ and $W$ preserving the curvature condition. We will view $S^{2k+1}$ as a product $[0,T] \times S^k_a \times S^k_b$ for some $T>>0$ (to be determined later), or more accurately as $$\Bigl(S^k_a \times \{0\}\Bigr)\cup \Bigl((0,T) \times S^k_a \times S^k_b\Bigr) \cup \Bigl(\{0\} \times S^k_b\Bigr).$$

The desired metric will be defined in three pieces: the `ends', which correspond to $t \in [0,T/9]$ and $t \in [8T/9,T]$; a middle region corresponding to $t \in [T/3,2T/3]$ in which the metric is chosen so as to facilitate surgery on $W$; and two `transition phases' corresponding to $t \in [T/9,T/3]$ and $t\in [2T/3,8T/9]$, in which the metric interpolates between those in the middle and end phases.  

For $t \in [0,T/9]$ we will choose the metric to be a warped product $$dt^2+\rho_0^2ds^2_k+q^2(t)ds^2_k.$$ Here $\rho_0$ is a small constant to be chosen later, and $q$ is a smooth function with the following properties: $q''(t) \le 0$; $q'(t) \in [0,1]$; $q(t)=\sin(t)$ if $t \in [0,\zeta],$ for some choice of $\zeta \in (0,1/10)$; $q(t)=\rho_1$, a small constant to be chosen later, for all $t \in [2\zeta,T/9]$. The form of $q$ close to $t=0$ ensures that the metric is smooth there.

For future reference, let us note the Ricci curvatures of the following single and double warped products. For a real interval $I$, the Ricci curvatures of the warped product $(I \times S^n;\ dt^2+p^2(t)ds^2_n)$ are given by
\begin{align*}
Ric(\partial_t,\partial_t)&=-np''/p, \\
Ric(U/p,U/p)&=-p''/p+(n-1)(1-p'^2)/p^2, 
\end{align*}
where $\partial_t$ is short for $\partial/\partial t,$ and where $U$ is a unit vector in the $S^n$-direction with respect to $ds^2_n$. The Ricci curvatures of the double warped product $$(I \times S^n \times S^m; dt^2+p^2(t)ds^2_n+q^2(t)ds^2_m)$$ are given by
\begin{align*}
Ric(\partial_t,\partial_t)&=-np''/p-mq''/q \\ 
Ric(U/p,U/p)&=-p''/p+(n-1)(1-p'^2)/p^2-mp'q'/pq \\ 
Ric(V/q,V/q)&=-q''/q+(m-1)(1-q'^2)/q^2-np'q'/pq \\
Ric(U/p,V/q)&=0,
\end{align*} 
where $\partial_t$ and $U$ are as above, and $V$ is a unit vector in the $S^m$-direction with respect to $ds^2_m$.

It is now easy to see that the `end' metric $dt^2+\rho_0^2ds^2_k+q^2(t)ds^2_k$ for $S^{2k+1}$ will have $Ric \ge 0$ and $Ric_2>0.$ Moreover we have $Ric>0$ for $t \in [0,\zeta].$ 

At the other end, we proceed similarly, introducing a metric $dt^2+p^2(t)ds^2_k+\rho_1^2 ds^2_k.$ We assume that $p''(t) \le 0$ and $p'(t) \in [-1,0]$ for all $t \in [8T/9,T].$ We further assume that for some $\zeta' \in (0,1/10),$ $p$ takes the form $$p(t)=\cos\bigl(\frac{\pi}{2}+t-T\bigr)$$ for $t \in [T-\zeta',T],$ and $p(t)=\rho_0$ for $t\in [8T/9,T-2\zeta'].$ Again, the resulting metric has $Ric \ge 0,$ $Ric_2>0,$ and for $t$ close to $T$ it is Ricci positive.

Notice that these warped product `end' metrics are essentially independent of $T$, in the sense that different values of $T$ (provided they are sufficiently large) affect neither the choice of $\zeta$ and $\zeta'$, nor the form of the scaling functions for $t \in [0,2\zeta]$ and $t \in [T-2\zeta',T]$, nor the constant values of the scaling functions at other values of $t$.
\medskip

We next select the metric on the middle region, corresponding to $t \in [T/3,2T/3].$
\medskip

The basic idea for the region $[T/3,2T/3] \times S^k_a \times S^k_b$ is to design our metric so as to reflect the topology of the embedded sphere $W$ and its normal bundle. We will first consider $W$ as a submanifold of $\{T/2\} \times S^k_a \times S^k_b$ (to simplify the notation we will write $W \subset S^k_a \times S^k_b$), construct a $W$-adapted submersion metric $g$ on $S^k_a \times S^k_b$, then extend trivially as a product to the whole neighbourhood.

At this point let us recall a result of Vilms (see \cite{Vi} or \cite[9.59]{Be}). Given a fibre bundle with Lie structure group, there is a submersion metric on the total space of the bundle with totally geodesic fibres which is uniquely determined by the following data: a metric on the base, a metric on the fibre invariant under the action of the structure group, and a principal connection on the associated principal bundle. The submersion metric is complete if the base and fibre metrics are both complete.

To obtain the desired metric on $S^k_a \times S^k_b$, we first consider creating a submersion metric on the disc bundle $D_\perp(W):=D_\perp(\Gamma(f))\subset S^k_a \times S^k_b$ (introduced before Corollary \ref{tubular}). We will do this in such a way that we automatically obtain a smooth submersion metric on the double of this bundle, i.e. on $S^k_a \times S^k_b$. 

With the above Vilms result in mind, we begin our construction by specifying a base metric: on $W$ choose the pull-back of the unit round metric via $\sigma,$ i.e. set $\check{g}=\sigma^*(ds^2_k)$, where $\sigma:W \to S^k_a$ is the diffeomorphism introduced after Corollary \ref{decomp}. 

For the connection, observe that our background metric $g_0$ induces Riemannian metrics in the fibres of $E(W):=E(\Gamma(f))$ (introduced after Corollary \ref{tubular}). This allows a reduction of the structure group to $\mathrm{SO}(k).$ Let $P$ denote the principal $\mathrm{SO}(k)$-bundle associated to $E(W)$. Notice that $P$ will also serve as the associated principal bundle for the disc bundle $D_\perp(W).$ We fix an arbitrary principal connection $\nabla$ on $P$. Now in the generic case where $d \neq 0$, $P$ must be non-trivial by Lemma \ref{Euler_number}, so consequently no choice of principal connection can be globally flat. The metric we will construct on $S^k_a \times S^k_b=D_\perp(W) \cup D_\perp(W)$ therefore cannot be a product metric in this case.

For the fibres we choose a rotationally symmetric metric $\hat{g}$ so as to be invariant under the structure group action. As such a metric is a single warped product, it follows from the above curvature formulas that the Ricci curvature will be positive provided the warping function is everywhere concave down. In fact we will choose this function so that the double approximates a long, thin, capped cylinder, which displays $\Z_2$-symmetry with respect to interchanging the poles. The precise form of this warping function will be crucial, and we will discuss this below.

We can now define a metric on $D_\perp(W)$ via the Vilms construction, starting with the data $\check g$, $\hat g$ and $(P,\nabla)$ above. For convenience, let us denote this metric by $$\nu(\check{g},\hat{g},\nabla),$$ and similarly for other metrics we will consider which arise from the Vilms construction. 

The desired metric $g$ on $S^k_a \times S^k_b$ is then the double of $\nu(\check{g},\hat{g},\nabla).$
\medskip

We now present the finer details of the fibre metric. Given that we want to work with the double of the disc bundle $D_\perp(W)$, it will be convenient to define the fibre metric on the whole of the sphere bundle fibre $S_b^k$, with the disc metric being `half' of this.
 
This standard fibre metric will take the form $ds^2+\psi^2(s)ds^2_{k-1}$ for $s \in [0,\Lambda],$ where the function $\psi$ is given by $$\psi(s)=\Lambda \theta(s/\Lambda)$$ for some large constant $\Lambda$ and a function $\theta$ defined below. Thus in terms of the notation introduced above we have $$g=\nu(\check g,ds^2+\psi^2(s)ds^2_{k-1},\nabla).$$ 
We will specify a value for $\Lambda$ in Section \ref{sec:surgery}, but assume for the moment that a suitable value for this constant has been fixed, and that $\Lambda>\pi.$ We therefore need to describe the function $\theta(x)$ where $x \in [0,1].$ We proceed as follows.

Consider the functions $$\frac{1}{\pi\Lambda}\sin(\pi\Lambda x)$$ for $x \in [0,\frac{1}{2\Lambda}],$ $$\frac{1}{\pi\Lambda}\sin\bigl(\pi\Lambda(x-1+\frac{1}{\Lambda})\bigr)$$ for $x \in [1-\frac{1}{2\Lambda},1],$ and $$c\sin(\pi(x+1)/3)$$ for $x \in [-1,2].$ Here $c>0$ is a very small constant. Notice that the first two functions taken together on the specified intervals display symmetry about $x=1/2,$ as does the third function. It is clear that for $c$ sufficiently small (depending on $\Lambda$), the graphs of the first and third functions above must intersect at a unique point $x_1 \in (0,1/2)$, and the graphs of the second and third must intersect at a point $x_2 \in (1/2,1).$ Moreover given any $\epsilon>0$, by choosing $c$ smaller if necessary, we can ensure that $x_1 \in (0,\epsilon)$ and $x_2 \in (1-\epsilon,1).$

Let $\theta_0(x)$ be the following continuous, piecewise-smooth function:
\[
\theta_0(x)=
\begin{cases}
\frac{1}{\pi\Lambda}\sin(\pi\Lambda x) & \text{if $x\in [0,x_1]$} \\
c\sin(\pi(x+1)/3) & \text{if $x \in (x_1,x_2)$} \\
\frac{1}{\pi\Lambda}\sin\bigl(\pi\Lambda(x-1+\frac{1}{\Lambda})\bigr) & \text{if $x\in [x_2,1].$} \\
\end{cases}
\]

In order to obtain $\theta$ we smooth $\theta_0.$ Specifically, we deform $\theta_0$ in very small neighbourhoods of $x=x_1$ and $x=x_2$ in such a way that the strict concavity of the function is preserved. It is clear that we can make such adjustments keeping $|\theta(x)-\theta_0(x)|$ arbitrarily small, and so that $\theta'(x)$ changes approximately linearly between its values on either side of the deformation regions. Moreover, to create $\Z_2$-symmetry about $x=1/2,$ we will assume that the two deformations performed will be mirror-images of one-another.

We remark that the `ends' of the metric $ds^2+\psi^2(s)ds^2_{k-1}$, corresponding to $s\approx 0$ and $s \approx \Lambda,$ are round of radius $1/\pi.$

With an eye towards curvature computations, we prove the following
\begin{lemma}\label{theta_est}
For any suitably small choice of $c>0,$ $\theta_0$ can be smoothed to a function $\theta(x)$ which satisfies
\begin{enumerate}
\item $-\theta''/\theta >1$;
\item $(1-\theta'^2)/\theta^2>1;$
\end{enumerate}
for all $x \in [0,1].$
\end{lemma}

\begin{proof}
We first check the above inequalities hold for $\theta_0$ on the intervals $[0,x_1)$ and $(x_1,x_2),$ then argue that these inequalities must continue to hold when $\theta_0$ is smoothed to $\theta.$ By the $\Z_2$-symmetry, this will then suffice to establish the Lemma.

For $x \in [0,x_1)$ we have $-\theta''_0/\theta_0=\pi^2\Lambda^2>1.$ For $x \in (x_1,x_2)$ we have $$-\theta''_0/\theta_0=\pi^2/9>1.$$

Given that $$\lim_{x \to x_1^-}\theta_0(x)=\cos(\pi\Lambda x_1)\,\, \text{ and }\lim_{x \to x_1^+}\theta_0(x)=\frac{\pi c}{3}\cos(\pi(x_1+1)/3),$$ we see that provided $c$ is sufficiently small,
the smoothing process must locally make the second derivative more negative, while keeping the function values arbitrarily close to the original $\theta_0$ values. Thus it is clear that we can perform the smoothing so that the inequality (1) holds.

For the second inequality we have $(1-\theta_0'^2)/\theta_0^2=\pi^2\Lambda^2$ when $x \in [0,x_1).$ For $x \in (x_1,x_2)$ we have $$\frac{1-\theta_0'^2(x)}{\theta_0^2(x)}=\frac{1-\frac{\pi^2}{9}c^2\cos^2(\pi(x+1)/3)}{c^2\sin^2(\pi(x+1)/3)},$$ and clearly this exceeds 1 if $c$ is small. 

As $\theta_0$ is smoothed to $\theta$, provided the interval of smoothing is sufficiently short, $(1-\theta'^2)/\theta^2$ will interpolate between its values on either side of the deformation to within any desired degree of accuracy. As the inequality holds on either side, it will therefore hold throughout.
\end{proof}

We immediately deduce
\begin{corollary}
\label{psi_estimates}
For any suitably small value of $c$, the following inequalities hold for all $s \in [0,\Lambda]$:
\begin{enumerate}
\item  $-\psi''/\psi >1/\Lambda^2$;
\item $(1-\psi'^2)/\psi^2>1/\Lambda^2.$
\end{enumerate}
\end{corollary}

Note that the required values for $c$ and $\Lambda$ needed to construct $\psi$ will emerge naturally from the surgery considerations in Section \ref{sec:surgery}.

In order to state the Ricci curvature formulas for the metric $g$ on $S^k_a \times S^k_b$, we introduce a Riemannian $S^{k-1}$-bundle over $W$ which will serve as a sort of `comparison bundle' to $\bar{\pi}:S^k_a \times S^k_b \to W.$ This bundle, which we will denote $\mathcal B$, is the $S^{k-1}$-bundle associated to the principal $\mathrm{SO}(k)$-bundle $P$, equipped with the submersion metric $g_B:=\nu(\check{g},ds^2_{k-1},\nabla).$ (Recall that this is the metric determined via the Vilms construction by the base metric $\check{g}$, the fibre metric $ds^2_{k-1},$ and principal connection $\nabla.$) The O'Neill formulas for the curvature of Riemannian submersions then describe the Ricci curvatures of $\mathcal B$ in terms of the corresponding $A$-tensor, see \cite[\S9C]{Be} for details. Note that $\mathcal B$ is essentially a sub-bundle of $\bar{\pi}:S^k_a \times S^k_b \to W$ corresponding to a constant value of $s \in (0,\Lambda)$, and up to fibre scaling, the metric $g_B$ agrees with the metric induced on $\mathcal B$ by $g$. In particular a vector $V$ tangent to a fibre of $\mathcal B$ can (and in the Lemma below will) also be viewed as tangent to a fibre of $\bar{\pi}.$ 

\begin{lemma}
\label{ambient_curvature}
The Ricci curvatures of the metric $g$ on $S^k_a \times S^k_b$ are as follows:
\begin{align*}
Ric(V/\psi,V/\psi)&=-\frac{\psi''}{\psi}+(k-2)\frac{1-\psi'^2}{\psi^2}+\psi^2(AV,AV); \\
Ric(\partial_s,\partial_s)&=-(k-1)\frac{\psi''}{\psi}; \\
Ric(X,X)&=(k-1)-2\psi^2(A_X,A_X); \\
Ric(V/\psi,X)&=-\psi(\check{\delta}A(X),V); \\
\end{align*}
with all other mixed Ricci curvature terms vanishing. Here $V$ is a unit vector for $(S^{k-1},ds^2_{k-1}),$ where $S^{k-1}$ is a cross-section of a fibre $S^k_b$ corresponding to any fixed $s\in(0,\Lambda),$ $\partial_s$ is shorthand for $\partial/\partial s$, and $X$ is a unit vector tangent to the base $W$. The $A$-tensor terms are those for the comparison bundle $({\mathcal B},g_B)$ (which in particular are independent of $\psi$).
\end{lemma}

For details behind these and many other curvature formulas appearing in this paper, see the Appendix.

Combining Lemma \ref{ambient_curvature} with Corollary \ref{psi_estimates}, the following result is evident:
\begin{proposition}\label{first_c_restriction}
There is a $c_0>0$ (depending on $\Lambda$ and $\nabla$) such that for each $c \in (0,c_0)$, the resulting function $\psi$ yields a metric $g$ which has all Ricci curvatures strictly positive.
\end{proposition}

We will assume from now on that our choices of $\Lambda$ and $c$ ensure that $Ric(g)>0.$ 

To complete the construction of middle-phase metric on $[T/3,2T/3] \times S^k_a \times S^k_b$ we simply form the product $dt^2+g.$
\bigskip

It remains therefore to describe the transition phase metrics corresponding to $t \in [T/9,T/3]$ and $t \in [2T/3,8T/9].$
\medskip

For the transition phases, and indeed also for surgery considerations in the next section, the following basic feature of $k$-positive Ricci curvature will be very useful:

\begin{lemma} 
\label{transition} 
Given a path of $Ric_k>0$ metrics $g_t$, $t \in [0,1]$, on a closed manifold $M$, there exists $\chi \ge 1$ such that the metric $dy^2+g_{L(y)}$ on $[0,\chi] \times M$ has $Ric_{k+1}>0,$ where $L$ denotes the linear function $[0,\chi] \to [0,1]$ given by $L(y)=y/\chi.$ In particular, any path of Ricci positive metrics on $M$ can be stretched in this way to give a $Ric_2>0$ metric on $[0,\chi] \times M.$ 
\end{lemma}

\begin{proof} 
We begin by observing that for fixed $t=t_1,$ the product metric $dt^2+g_{t_1}$ has $Ric_{k+1}>0$: $Ric_k>0$ restricted to $M$-directions, and zero Ricci curvature in the $t$-direction.

Given $\epsilon>0,$ by choosing the value of $\chi$ sufficiently large, we can ensure that the metric $dy^2+g_{L(y)}$ in a neighbourhood of $y=\chi t_1$ is $C^2$ $\epsilon$-close to the $Ric_{k+1}>0$ product metric $dy^2+g_{t_1}.$ As curvature only depends on the metric in a $C^2$-sense, by the openness of the $Ric_{k+1}>0$ condition, we conclude that the metric $dy^2+g_{L(y)}$ also has $Ric_{k+1}>0$ in this neighbourhood for $\epsilon$ sufficiently small. By the compactness of the interval $[0,\chi]$, we can therefore choose a value for $\chi$ which gives $Ric_{k+1}>0$ for the metric $dy^2+g_{L(y)}$ globally. 
\end{proof}

By Lemma \ref{transition}, it suffices to show that there is a path of Ricci positive metrics on $S^k_a \times S^k_b$ linking $g$ above with the product metric $\rho_0^2 ds^2_k+\rho_1^2 ds^2_k$ on $S^k_a \times S^k_b$ for suitable $\rho_0,\rho_1>0.$ The transition phase breaks into four sub-phases: a re-scaling phase, a twisting phase, a shape-changing phase, and then a further re-scaling phase.
\bigskip

\noindent $\bullet$ {\it The first rescaling.}
\medskip

\begin{lemma}
\label{first_rescale}
For any $\rho_0,\rho_1>0$ there exists a Ricci positive path of metrics joining $\rho_0^2 ds^2_k+\rho_1^2 ds^2_k$ to $ds^2_k+\rho_1^2 ds^2_k.$
\end{lemma}

\begin{proof}
As the metric $\mu^2 ds^2_k+\rho_1^2 ds^2_k$ has positive Ricci curvature for all $\mu>0,$ by varying $\mu$ between $\rho_0$ and 1 we obtain a path of Ricci positive metrics linking $\rho_0^2 ds^2_k+\rho_1^2 ds^2_k$ to $ds^2_k+\rho_1^2 ds^2_k.$
\end{proof}

The metric in the above Lemma constitutes the first part of the transition phase.
\bigskip

\noindent $\bullet$ {\it The twisting phase.}
\medskip

The basic background result here is:
\begin{lemma}\label{can_var}
Consider a Riemannian submersion with totally geodesic fibres determined via the Vilms construction by a base metric, a fibre metric, and a connection on the associated principal bundle. If the base and fibre metrics are Ricci positive, then for any choice of connection there is a $\mu_0>0$ such that for any $\mu\in (0,\mu_0)$, rescaling the fibre metrics by a factor of $\mu^2$ yields a Ricci positive metric on the total space.
\end{lemma}

\begin{proof}
This is an immediate consequence of the canonical variation formulas, as presented in \cite[9.70]{Be}.
\end{proof}

Recall that the Ricci positive metric $g$ on $S^k_a \times S^k_b$ constructed for the middle region is a submersion metric with totally geodesic fibres, determined by the Ricci positive base metric $(W,\check g)$, a fibre metric $(S^k_b,ds^2+\psi^2(s)ds^2_{k-1}),$ and a principal connection $\nabla$ on the associated principal $SO(k)$-bundle $P$.

Notice that we can also view $g$ as a Riemannian submersion over $(S^k_a,ds^2_k)$ given by the projection map $\pi=\sigma \circ \bar{\pi}.$ This works because the original base metric on $W$ is the pull-back of $ds^2_k$ via $\sigma.$ We are therefore justified in writing $$g=\nu(ds^2_k,ds^2+\psi^2(s)ds^2_{k-1},\nabla),$$ with the projection map $\pi$ assumed. This shift in our viewpoint will aid comparison with product metrics on $S^k_a \times S^k_b$, which of course are submersion metrics for $\pi$ asssociated to the trivial (flat) principal $SO(k+1)$-connection $\nabla_{triv}.$

Joining $g$ to product metrics on $S^k_a \times S^k_b$ via a path of Ricci positive submersion metrics will, in particular, require an `untwisting' of $\nabla.$ However this is not possible as things currently stand, as the associated principal $SO(k)$-bundle $P$ is typically non-trivial.

In order to address this, we first have to enlarge the isometry group of the fibres from $SO(k) \times \Z_2$ to $SO(k+1),$ by replacing the metric $ds^2+\psi^2(s)ds^2_{k-1}$ by the round metric $ds^2_k.$ Thus we will consider making a suitable deformation to the scaling function $\psi$. 

If the unit round metric is used in the Vilms construction with the same base metric and principal connection, we obtain a submersion metric with round totally geodesic fibres. Call this resulting metric $g_1.$ As a consequence of the round fibres, for this new metric we now have an associated principal $SO(k+1)$-bundle $P_1,$ and a principal $SO(k+1)$-connection $\nabla_1.$ We can therefore write $g_1=\nu(ds^2_k,ds^2_k,\nabla_1).$

It is easy to see that $P_1$ is a trivial bundle, since it is associated to the trivial $S^k$-bundle $\pi=\sigma \circ \bar{\pi}$ (whereas $P$ is associated to the generally non-trivial $S^{k-1}$-`bundle of equators' within $S^k_a \times S^k_b$ determined by $W$). Hence we have two principal $SO(k+1)$-connections, $\nabla_1$ and $\nabla_{triv}$ on the trivial $SO(k+1)$-bundle over $S^k_a.$

As the space of principal connections on a given principal bundle is a contractible space, we can find a path of principal connections joining $\nabla_{triv}$ to $\nabla_1.$ This path, together with the fixed base $(S^k_a,ds^2_k)$ and fibre $(S^k_b,ds^2_k)$ produces a path of submersion metrics via the Vilms construction with round totally geodesic fibres linking $g_0$ to $g_1.$

\begin{lemma}
\label{iota_1}
There exists $\iota_1\in(0,1)$ such that for any $\rho_1 \in (0,\iota_1)$, rescaling the fibres of each submersion metric in the above path linking $g_0$ to $g_1$ by the factor $\rho_1^2$ yields a path of Ricci positive metrics.
\end{lemma}

\begin{proof} 
This follows from Lemma \ref{can_var} together with the compactness of the path domain.
\end{proof}

Let us denote the metric obtained from $g$ by re-scaling its fibres by a factor $\mu^2$ by $g(\mu),$ and similarly for other submersion metrics, so Lemma \ref{iota_1} guarantees a Ricci positive path from $g_0(\rho_1)$ to $g_1(\rho_1).$
\bigskip

\noindent $\bullet$ {\it The shape changing phase}
\medskip

We now consider changing the shape of the fibres from $\rho_1^2 ds^2_k$ to $\rho_1^2(ds^2+\psi^2(s)ds^2_{k-1}).$ 
To this end, we will need to describe a path from the round metric to $ds^2+\psi^2(s)ds^2_{k-1}.$ Now the latter metric is a warped product over the interval $[0,\Lambda].$ In order to describe the path of metrics it will be convenient to express the round metric similarly as a warped product over $[0,\Lambda]$. The natural way to do this is as follows: $$ds^2+\frac{\Lambda^2}{\pi^2}\sin^2\Bigl(\frac{\pi s}{\Lambda}\Bigr)ds^2_{k-1}.$$ Of course this is a round metric of radius $\Lambda/\pi$, i.e. $\Lambda^2\pi^{-2}ds^2_k.$ Given this representation of the round metric, we will now describe the desired path. Of course we must make sure that all the individual metrics along this path are Ricci positive. 

Let $\beta:[0,1] \to [0,1]$ be a smooth `bump' function satisfying $\beta(x)=0$ for $x \in [0,\epsilon],$ where $0<\epsilon <<1,$ $\beta(x)=1$ for $x \in [1-\epsilon,1],$ and $\beta'(x) \ge 0$ for all $x \in [0,1].$ Consider the function $H:[0,\Lambda] \times [0,1] \to \R_{\ge 0}$ given by $$H(s,x)=(1-\beta(x))\Bigl( \frac{\Lambda}{\pi}\sin\bigl(\frac{\pi s}{\Lambda}\bigr)\Bigr)+\beta(x)\psi(s).$$

\begin{lemma}\label{H}
For each $x \in [0,1],$ the warped product metric $ds^2+H^2(s,x)ds^2_{k-1}$ on $S^k$ is a smooth Ricci positive metric.
\end{lemma}

Clearly, this metric path begins (for $x \approx 0$) with a round metric of radius $\Lambda/\pi,$ and ends (for $x \approx 1$) with the metric $ds^2+\psi^2(s)ds^2_{k-1}.$ 

\begin{proof} 
Differentiating, we have 
\begin{align*}
\frac{\partial H}{\partial s}&=(1-\beta(x))\cos\bigl(\frac{\pi s}{\Lambda}\bigr)+\beta(x)\psi'(s), \\
\frac{\partial^2 H}{\partial s^2}&=(1-\beta(x))\Bigl(-\frac{\pi}{\Lambda}\sin\bigl(\frac{\pi s}{\Lambda}\bigr)\Bigr)+\beta(x)\psi''(s). \\
\end{align*}

To prove smoothness, by $\Z_2$-symmetry about $s=\Lambda/2$ it suffices to study $s \approx 0.$ Here we require oddness at $s=0$ and $\frac{\partial H}{\partial s}(0,x)=1$ for all $x$. Now a convex combination of odd functions is clearly odd, and $(\Lambda/\pi)\sin(\pi s/\Lambda)$ and $\psi(s)$ are individually odd. Moreover $$\frac{\partial H}{\partial s}(0,x)=(1-\beta(x)).1+\beta(x).1=1,$$ as required.

For Ricci positivity we observe that $H_{ss}/H<0$ for $s \ge 0$ (using l'H\^opitals's rule to examine the case $s=0$). As $(1-H_s^2)/H^2$ is non-negative (again, using l'H\^opital to investigate $s=0$), Ricci positivity now follows from the single warped product formulas at the start of this section.
\end{proof}

\begin{lemma}\label{iota_2}
There exists $\iota_2>0$ such that for any $\rho_2\in (0,\iota_2),$ the path of totally geodesic submersion metrics (parametrized by $x \in [0,1]$) with fibre $$(S^k_b, \rho^2_2(ds^2+H^2(s,x))ds^2_{k-1}),$$ base $(S^k_a,ds^2_k)$ and associated principal connection $\nabla$, is a path of Ricci positive metrics. In other words we have a path of Ricci positive metrics joining $g_1(\Lambda\rho_2/\pi)$ to $g(\rho_2).$
\end{lemma}

\begin{proof} 
This is evident by combining Lemma \ref{H} with Lemma \ref{can_var}, bearing in mind that $ds^2+H^2(s,0)ds^2_{k-1}$ is a round metric of radius $\Lambda/\pi$.
\end{proof}

We now need to factor $\rho_2$ in the above Lemma into our choice of $\rho_1,$ so that this $\rho_1$ will work for both Lemmas 
\ref{iota_1} and \ref{iota_2}.

\begin{corollary}
\label{path_1}
For any $\rho_1 \in (0,\min\{\iota_1,\iota_2\Lambda/\pi\}),$ there is a Ricci positive path of totally geodesic submersion metrics linking $g_0(\rho_1)=ds^2_k+\rho_1^2 ds^2_k$ and $g(\pi\rho_1/\Lambda).$
\end{corollary}
\medskip

\noindent $\bullet$ {\it Second scaling phase.}
\medskip

\begin{lemma}
\label{path_2}
There exists a path of Ricci positive metrics linking $g(\pi\rho_1/\Lambda)$ and $g$.
\end{lemma}

\begin{proof}
Consider the path of submersion metrics taking the form $g(\mu)$ for $\mu \in (0,1].$ Since $Ric(g)>0,$ we obtain Ricci positivity for each metric of this form from Lemma \ref{can_var}. The result now follows from the fact that $\pi\rho_1/\Lambda<1.$ This inequality is an elementary consequence of the fact that $\rho_1<1$ (see Lemma \ref{iota_1}) and $\Lambda>\pi.$ 
\end{proof}

\begin{proposition}
There exists a Ricci positive path of metrics joining $\rho_0^2 ds^2_k+\rho_1^2 ds^2_k$ to the metric $g$.
\end{proposition}

\begin{proof}
Simply concatenate the path in Lemma \ref{first_rescale} with those described by Corollary \ref{path_1} and Lemma \ref{path_2}.
\end{proof} 

Using Lemma \ref{transition} we obtain

\begin{corollary}\label{T_and_rho_1}
If $T$ is chosen sufficiently large, and then $\rho_1$ chosen sufficiently small, there exists a smooth $Ric_2>0$ metric on $(0,2T/3] \times S^k_a \times S^k_b$ which agrees with the metric $dt^2+\rho_0^2 ds^2_k+q^2(t)ds^2_k$ for $t \in (0,T/9],$ and with $dt^2+g$ for $t \in [T/3,2T/3].$
\end{corollary}

It remains to construct the part of the $Ric_2>0$ metric on $S^{2k+1}$ corresponding to $t \in [2T/3,8T/9],$ to join smoothly with $dt^2+g$ for smaller $t$, and the metric $dt^2+p^2(t)ds^2_k+\rho_1^2 ds^2_k$ at larger $t$. However given the construction above for the transition phase $t \in [T/9,T/3],$ it is easy to see that the desired metric is essentially the mirror image of the above.


\bigskip\medskip
\section{Surgery on links}
\label{sec:surgery}
\bigskip\medskip

We begin by performing Ricci positive surgery on the sphere $Z:=S^k_a \times \{0\}.$ Recall that for $t \in [0,\zeta],$ the ambient metric takes the form $dt^2+\rho_0^2 ds^2_k+\sin^2(t) ds^2_k.$ The following result essentially follows from \cite[Lemma 1]{SY}. 
\begin{proposition}
\label{V_surgery}
If $\rho_0$ is sufficiently small (given $\zeta$), then surgery can be performed on $Z$ within the region corresponding to $t \in [0,\zeta/2]$, preserving Ricci positivity.
\end{proposition}

\begin{remark}
The surgery we are dealing with here is the `obvious' surgery, in the sense that the normal bundle of $Z$ has an evident product structure, which we then use to complete the surgery. However, given a different trivialization of the normal bundle, one can equally perform surgery using that datum. Topologically this can result in a different manifold, and we must consider such possibilities for each surgery we perform. For the moment we will not worry about this point: we will work with the obvious surgeries in the short term, and address the situation for alternative normal bundle trivializations at the end of the section.
\end{remark}

Given that surgery on $Z$ is straightforward, the main focus of this section will therefore be concerned with performing $Ric_2>0$ surgery on the sphere $W$. Recall that $W$ is embedded in $\{T/2\} \times S^k_a \times S^k_b,$ and that the metric on $[T/3,2T/3] \times S^k_a \times S^k_b$ is $dt^2+g.$

Our first task is to identify the normal discs which will be part of the surgery procedure. It is clear that at any point $w=(T/2,x,f(x)) \in W$, the submanifold $[T/3,2T/3] \times \{x\} \times S^k_b$ meets $W$ orthogonally at $w$. Moreover, as the fibres of the submersion metric $g$ are totally geodesic, the exponential map defined by $g$ in $\{T/2\} \times S^k_a \times S^k_b$ gives a local diffeomorphism of the normal space to $W\subset \{T/2\} \times S^k_a \times S^k_b$ at $w$ onto the fibre $\{T/2\} \times \{x\} \times S^k_b.$ Thus the desired normal disc to $W$ at $w$ within the ambient manifold will be some choice of disc $D^{k+1} \subset [T/3,2T/3] \times \{x\} \times S^k_b$.

Consider the normal distance sphere $S^k$ about $w \in W$ at some radius $R$. In order for this sphere to lie wholly within the region $[T/3,2T/3] \times \{x\} \times S^k_b$, we need to place a lower bound on $T$, namely 
\begin{equation}
T>6R.
\end{equation}
Recalling that we have a longitudinal parameter $s$ on each $S^k_b$, observe that the normal distance sphere corresponds to a semi-circle in the $(t,s)$-plane. Specifically, let $\gamma(r),$ $r \in [0,\pi R]$, denote the radius $R$ semi-circle $$\gamma(r)=\bigl(R\cos(r/R)+\frac{T}{2},R\sin(r/R)\bigr).$$ Notice that each point of this curve for $r \in (0,\pi R)$ picks out a copy of $S^{k-1}$ in the appropriate fibre $S^k_b$ corresponding to the $s$-parameter value $R\sin(r/R)$, i.e. $$S^{k-1} \subset \Bigl\{R\cos(r/R)+\frac{T}{2}\Bigr\} \times \{x\} \times S^k_b \subset [T/3,2T/3] \times S^k_a \times S^k_b.$$ At $r=0,\pi R$ we obtain single points instead, so topologically, the curve $\gamma$ determines a $k$-sphere as claimed. As a distance sphere with respect to a smooth metric, this is clearly a smooth submanifold of $S^{2k+1}.$

It turns out that $R$ will need to be sufficiently large in order to perform $Ric_2>0$ surgery on $W$. Clearly there has to be a relationship between $R$ and $\Lambda.$ It will be convenient to set
\begin{equation}
\Lambda=\frac{5R}{2}.
\end{equation}

Collectively, i.e. considering all possible $(x,f(x)) \in W$ (or equivalently all $x \in S^k_a$), these normal distance spheres form a trivial $S^k$-bundle over $W$, since the normal (vector) bundle to $W$ in $S^{2k+1}$ is trivial by Lemma \ref{s-trivial}. Metrically, however, this $S^k$-bundle is generally not trivial (that is, not a metric product), as the $S^{k-1}$ cross-sections of the $S^k$ fibres (corresponding to constant values of $r$) are subject to twisting as one moves around the base $W$, as determined by the principal $SO(k)$-connection $\nabla$ associated to $g$.

For surgery purposes we want to consider not just the copies of the normal sphere $S^k$ (for a given point $(x,f(x)) \in W$), but the family of `concentric' normal spheres which (inwardly) foliate the normal disc bounded by this $S^k.$

For $\tau\ge 0$ we therefore consider the inward $\tau$-equidistant semi-circle $$\gamma(r,\tau)=\bigl( (R-\tau)\cos(r/R)+\frac{T}{2},(R-\tau)\sin(r/R)\bigr)$$ for $r \in [0,\pi R].$ It is easily checked that the induced metric on the corresponding sphere $S^k_\tau$ is $$\Bigl(\frac{R-\tau}{R}\Bigr)^2dr^2+\psi^2\bigl((R-\tau)\sin(r/R)\bigr)ds^2_{k-1}.$$ (Setting $u=r(R-\tau)/R$ and re-expressing the metric in terms of $u$, we again see that the metric is smooth for all $\tau \in [0,R).$)

Notice now that $(r,\tau)$ together form a coordinate system on part of the $(t,s)$-plane - the part relevant to our surgery considerations. It will make sense to switch to using $(r,\tau)$ coordinates from now on.

Consider the normal $D^{k+1}$ disc bundle about $W$ where the fibre-disc boundaries are determined by $\gamma(r)=\gamma(r,0)$. 
Notice that $P$, the principal $SO(k)$-bundle associated to the normal bundle of $W \subset S^k_a \times S^k_b$, is also an associated principal bundle to our normal $D^{k+1}$-bundle over $W$, (with $SO(k)$ acting trivially in the $t$-direction). In particular, $\nabla$ is the appropriate principal connection for the normal $D^{k+1}$-bundle.

The ambient metric in this disc bundle can then be described as a submersion metric over $(W,\sigma^*(ds^2_k))$, with principal connection $\nabla$, and with fibre metric given in our new $(r,\tau)$ coordinate system by $$d\tau^2+\Bigl(\frac{R-\tau}{R}\Bigr)^2dr^2+\psi^2\bigl((R-\tau)\sin(r/R)\bigr) ds^2_{k-1}.$$

We will initially consider $\tau \in [0,1].$ Assuming $R>>1,$ this will not cause us any problems.

Our first task in the surgery process is to `slow down' the evolution of concentric normal sphere metrics as $\tau$ increases, to a point where the normal sphere metrics are independent of $\tau$. We achieve this by introducing a function $\phi:[0,1]\to [0,1]$ given by $\phi(\tau)=\tau$ for $t\in [0,1/10],$ $\phi(\tau)=1/2$ for $\tau \in [9/10,1],$ with $\phi' \in [0,1]$ and $\phi'' \in (-3/2,0].$ (It is clear that such a function exists.) Our strategy will be to replace the metric on the normal sphere $S^k_\tau$ by that for $S^k_{\phi(\tau)},$ so for $\tau \in [9/10,1]$ the metric does not vary with $\tau$. Thus we want to consider instead the normal metric
\begin{equation}\label{normal_sphere_metric}
d\tau^2+\Bigl(\frac{R-\phi(\tau)}{R}\Bigr)^2 dr^2+\bar{\psi}^2(r,\tau)ds^2_{k-1},
\end{equation}
where we have set $$\bar{\psi}(r, \tau):=\psi\bigl((R-\phi(\tau))\sin(r/R)\bigr).$$ 
Since $\tau \in [0,1]$ here, we are actually dealing with a normal annulus bundle, which is a subbundle of the normal disc bundle.

\begin{remark}
Setting $\Lambda=5R/2$ means that $(R-\tau)\sin(r/R)<\Lambda/2,$ which ensures the concavity of $\bar{\psi}.$ This will be important for curvature considerations.
\end{remark}

To allow greater flexibility, we will also allow base scalings: instead of a fixed metric $\sigma^*(ds^2_k)$ on $W$, we will instead consider $f^2(\tau)\sigma^*(ds^2_k)$ for some function $f$. To join smoothly with the ambient metric outside the normal disc bundle, we clearly need $f(\tau)=1$ for $\tau \approx 0.$

Let us denote by $g_\phi$ the metric on the normal annulus bundle to $W$ corresponding to $\tau \in [0,1],$ determined by the fibre metric (\ref{normal_sphere_metric}), base metric $f^2(\tau)\sigma^*(ds^2_k),$ and the principal connection $\nabla$. Notice that we can view this as a submersion metric over $(W,\sigma^*(ds^2_k))$ followed by a base direction rescaling. A calculation described in the Appendix shows the following.

\begin{proposition}
\label{ric(g_phi)}
Suppose that the vector $V$ is as in Lemma \ref{ambient_curvature}, and that $\{X_i\}$ is an orthonormal frame in $(W,\sigma^*(ds^2_k))$. Then the metric $g_\phi$ has Ricci curvatures given by
\begin{align*}
Ric\Bigl(\frac{1}{1-R^{-1}\phi}\frac{\partial}{\partial r},&\frac{1}{1-R^{-1}\phi}\frac{\partial}{\partial r}\Bigr)=\frac{\phi''}{R-\phi}+k\frac{\phi'}{R-\phi}\frac{f'}{f}-\frac{k-1}{(1-R^{-1}\phi)^2}\frac{\bar{\psi}_{rr}}{\bar{\psi}} \\
&\hskip 3.5cm +(k-1)\frac{\phi'}{R-\phi}\frac{\bar{\psi}_\tau}{\bar{\psi}},\\
Ric(V/\bar\psi,V/\bar\psi)=&\frac{k-2}{\bar{\psi}^2}\Bigl(1-\bar{\psi}^2_\tau-\frac{\bar{\psi}^2_r}{(1-R^{-1}\phi)^2}\Bigr)-\frac{1}{\bar\psi}\Bigl(\bar{\psi}_{\tau\tau}+\frac{\bar{\psi}_{rr}}{(1-R^{-1}\phi)^2}\Bigr) \\
& +\Bigl(\frac{\phi'}{R-\phi}-k\frac{f'}{f}\Bigr)\frac{\bar{\psi}_\tau}{\bar\psi}+\frac{\bar{\psi}^2}{f^4}(AV,AV), \\
Ric(\partial/\partial \tau,\partial/\partial \tau)=&\frac{\phi''}{R-\phi}-k\frac{f''}{f}-(k-1)\frac{\bar{\psi}_{\tau\tau}}{\bar{\psi}},\\
Ric(X_i/f,X_i/f)=&\frac{k-1}{f^2}(1-f'^2)-\frac{f''}{f}+\frac{\phi'}{R-\phi}\frac{f'}{f}-(k-1)\frac{f'\bar{\psi}_\tau}{f\bar{\psi}}-2\frac{\bar{\psi}^2}{f^4}(A_{X_i},A_{X_i}),\\
Ric(X_i/f,V/\bar{\psi})=&-\frac{\bar\psi}{f^3}(\check{\delta}A(X_i),V).\\
Ric(X_i/f,X_j/f)=&-2\frac{\bar{\psi}^2}{f^4}(A_{X_i},A_{X_j}) \,\,\mathrm{ for }\, i \neq j,\\
Ric(V/\bar{\psi},V'/\bar{\psi})=&\frac{\bar{\psi}^2}{f^4}(AV,AV')\,\, \mathrm{if}\, ds^2_{k-1}(V,V')=0. \\
\end{align*}
All other mixed Ricci curvature terms vanish. The $A$-tensor terms are those for the comparison bundle $({\mathcal B},g_B).$
\end{proposition}

Let us now give more detail about the function $f(\tau)$ for $\tau \in [0,1].$ As noted above, we need $f(\tau)=1$ for $\tau$ close to zero. Suppose we are given a very small constant $\delta>0.$ (This constant will have to be chosen sufficiently small for later arguments to be valid, but this will not concern us for the moment.) Let us choose and fix a function $f$ such that $f''(\tau)\le 0$, $f''(\tau)=-\delta$ for $\tau \in [1/10,9/10],$ $f''(\tau)=0$ for $\tau \approx 1$, and $$\frac{f'}{f}>-\delta.$$ Over the interval $[1/10,9/10]$, $f'$ decreases by $4\delta/5$. It is clear that we can choose $f$ for the complementary values of $\tau$ in such a way that $f'(1)>-9\delta/10.$ In turn this means that $$f(1)>1-\frac{9\delta}{10}.$$ Therefore the minimum value of $f'/f$ for $\tau \in [0,1]$ occurs at $\tau=1$, and we have 
\begin{align*}
\frac{f'}{f}(1)&>-\frac{9\delta/10}{1-(9\delta/10)} \\
&=-\frac{9\delta}{10-9\delta} \\
&>-\delta,\,\,\, \text{ provided } \delta<1/9. \\
\end{align*}
Clearly we also have $-f''/f>\delta$ for $\tau \in [1/10,9/10],$ and $-f''/f \ge 0$ for all $\tau \in [0,1].$
\bigskip\medskip

\centerline{\it The curvature of $g_\phi$}
\bigskip

The major task at hand is to analyse the Ricci curvature formulas for $g_\phi$ in the Proposition above. To this end, we recall that the function $\psi$ (on which $\bar\psi$ is based) is defined, via a $C^0$ function $\theta_0$, in three pieces: see Section \ref{sec:metrics} for the precise details. As our analysis will in particular involve estimates for $\bar\psi$ and its derivatives, we will split our considerations in pieces corresponding to those for $\theta_0.$ Of course we must also take into account the smoothing of $\theta_0$ to $\theta$ needed to define $\psi,$ and ultimately $\bar\psi$.
\bigskip

\noindent $\bullet$ {\it Near the fibre ends.}
\medskip

Let us recall that the non-smooth points for $\theta_0(x)$ occur at $x_1 \in (0,1/2)$ and $x_2 \in (1/2,1).$ Moreover, we asserted in Section \ref{sec:metrics} that given any $\epsilon>0$ we can choose $c$ sufficiently small so that $x_1 \in (0,\epsilon)$ and $x_2 \in (1-\epsilon,1).$ For the moment we will ignore all smoothing issues.

For small $x\ge 0$ we have $\theta(x)=\theta_0(x)=\pi^{-1}\Lambda^{-1}\sin(\pi\Lambda x).$ So for suitably small $r\ge 0$, where `suitably small' here is dependent (via the defining expression for $\bar{\psi}$) on $\tau \in [0,1]$, we have $$\bar{\psi}(r,\tau)=\frac{1}{\pi}\sin\bigl(\pi(R-\phi)\sin(r/R)\bigr),$$ and thus the following derivative formulas apply: 
\begin{align*}
-\frac{1}{(1-R^{-1}\phi)^2}\frac{\bar{\psi}_{rr}}{\bar\psi}&=\pi^2\cos^2(r/R)+\frac{\pi}{R-\phi}\sin(r/R)\cot\bigl(\pi(R-\phi)\sin(r/R)\bigr),\\
-\frac{\bar{\psi}_{\tau\tau}}{\bar\psi}&=\pi\phi''\sin(r/R)\cot\bigl(\pi(R-\phi)\sin(r/R)\bigr)+\pi^2\phi'^2\sin^2(r/R),\\
\frac{\bar{\psi}_\tau}{\bar\psi}&=-\pi\phi'\sin(r/R)\cot\bigl(\pi(R-\phi)\sin(r/R)\bigr),\\ \\
\frac{1-\bar{\psi}^2_\tau-(1-R^{-1}\phi)^{-2}\bar{\psi}^2_r}{\bar{\psi}^2}&=\pi^2\frac{1-\cos^2(\pi(R-\phi)\sin(r/R))(\cos^2(r/R)+\phi'^2\sin^2(r/R))}{\sin^2(\pi(R-\phi)\sin(r/R))}.\\
\end{align*}

Let us analyse the above expressions. 
\medskip

The final expression is bounded below by $\pi^2$, which can be seen by replacing $\phi'^2$ by 1.

The the first term in the expression for $-(1-R^{-1}\phi)^{-2}\bar{\psi}_{rr}/\bar{\psi}$ is clearly positive provided $r<\pi R/2.$ The behaviour of the second term in this expression as $r \to 0$ is governed by $$\lim_{r \to 0} \sin(r/R)\cot(\pi(R-\phi)\sin(r/R)).$$ By l'H\^opital's rule, we easily compute $$\lim_{r \to 0} \sin(r/R)\cot\bigl(\pi(R-\phi)\sin(r/R)\bigr)=\frac{1}{\pi(R-\phi)}.$$ Therefore overall, $$-\frac{1}{(1-R^{-1}\phi)^2}\lim_{r \to 0}\frac{\bar{\psi}_{rr}}{\bar\psi}=\pi^2+\frac{1}{(R-\phi)^2}.$$ Thus for $r$ sufficiently small, the expression for $-(1-R^{-1}\phi)^{-1}\bar{\psi}_{rr}/\bar{\psi}$ is positive. 

For the second expression we have, using the above, $$\lim_{r \to 0}\Bigl(-\frac{\bar{\psi}_{\tau\tau}}{\bar\psi}\Bigr)=\pi\phi''\frac{1}{\pi(R-\phi)}+0=\frac{\phi''}{R-\phi}.$$ This is therefore non-positive. For $\tau\in [0,1/10]\cup[9/10,1]$ we have $\phi''=0,$ but for the in-between values of $\tau$ this is a bad term which has to be managed. Note however that given the choice of $\phi$ (and in particular the fact that $\phi''>-3/2$), this term is at worst $-2/R,$ (assuming $R$ is large).

Finally, for the remaining expression we have $$\lim_{r \to 0} \frac{\bar{\psi}_\tau}{\bar\psi}=-\pi\phi'\frac{1}{\pi(R-\phi)}=-\frac{\phi'}{R-\phi}.$$ Notice that this is at worst $-1/(R-\phi),$ and so we can again use a lower bound of $-2/R.$

To conclude this analysis, we note note that the limiting values computed above serve as good estimates for the relevant terms provided $\epsilon,$ or equivalently $c$, is chosen small enough, (as these constants control $x_1$, and in turn the values of $r$ for which the above derivative expressions are valid). Assuming this, we now use these estimates to obtain lower bounds on the Ricci curvature terms.

We begin with $Ric(V/\bar{\psi},V/\bar{\psi})$. Using the crude estimates $\phi''>-2,$ $\phi'\le 1,$ $\phi<1,$ and $\delta<1,$ we obtain
\begin{align*}
Ric(V/\bar{\psi},V/\bar{\psi})\approx &(k-2)\pi^2+\frac{\phi''}{R-\phi}+\pi^2+\frac{1}{(R-\phi)^2}-\Bigl(\frac{\phi'}{R-\phi}\Bigr)^2 \\
&+k\frac{\phi'}{R-\phi}\frac{f'}{f}+\frac{\bar{\psi}^2}{f^4}(AV,AV),\\
\ge &(k-1)\pi^2-\frac{2}{R-\phi}+\frac{1}{(R-\phi)^2}-\frac{1}{(R-\phi)^2}-\frac{k\delta}{R-\phi},\\
\ge &(k-1)\pi^2-\frac{k+2}{R-1}.\\
\end{align*}
It is then easily checked that we obtain a positive lower bound (assuming that $k \ge 3$) if $R>2.$

In the $\tau$-direction we have
\begin{align*}
Ric(\partial/\partial \tau,\partial/\partial \tau)&\approx \frac{\phi''}{R-\phi}-k\frac{f''}{f}+(k-1)\frac{\phi''}{R-\phi},\\
&=k\Bigl(\frac{\phi''}{R-\phi}-\frac{f''}{f}\Bigr).\\
\end{align*}
Recalling our choices for $f$ and $\phi$ we have for $\tau \in [1/10,9/10]$
$$Ric(\partial/\partial \tau,\partial/\partial \tau)\ge k\Bigl(\delta-\frac{2}{R-1}\Bigr).$$
Thus this Ricci curvature will be positive at these values of $\tau$ if $\delta>2/(R-1),$ or equivalently if 
\begin{equation}
R>1+(2/\delta).
\end{equation}
For $\tau \in [0,1/10] \cup[9/10,1]$ we clearly have $Ric(\partial/\partial \tau,\partial/\partial \tau)\ge 0.$

In the $r$-direction we have
\begin{align*}
Ric\Bigl(\frac{1}{1-R^{-1}\phi}\frac{\partial}{\partial r},\frac{1}{1-R^{-1}\phi}\frac{\partial}{\partial r}\Bigr)\approx &\frac{\phi''}{R-\phi}+k\frac{\phi'}{R-\phi}\frac{f'}{f}+(k-1)\pi^2 \\
&+\frac{k-1}{(R-\phi)^2}-(k-1)\Bigl(\frac{\phi'}{R-\phi}\Bigr)^2,\\
\ge& -\frac{2}{R-1}+(k-1)\pi^2-k\frac{\delta}{R-1}.\\
\end{align*}
This will be positive provided $$(k-1)\pi^2>\frac{2+k}{R-1}.$$ Since $k \ge 3$ this will be satisfied if $R>2.$

Finally, we consider
\begin{align*}
Ric(X_i/f,X_i/f)&\approx (k-1)\frac{(1-f'^2)}{f^2}-\frac{f''}{f}+k\frac{\phi'}{R-\phi}\frac{f'}{f}-2\frac{\bar{\psi}^2}{f^4}(A_{X_i},A_{X_i}), \\
&\ge (k-1)(1-\delta^2)-\frac{k\delta}{R-1}-2\frac{\Lambda^2 c^2}{(1-\delta)^4}(A_{X_i},A_{X_i}).\\
\end{align*}
Clearly we can render the last term arbitrarily small by a suitable choice of $c$ (given the $A$-tensor term and a value for $\Lambda$), hence for positivity it suffices to show that $$(k-1)(1-\delta^2)>\frac{k\delta}{R-1}.$$ Since $\delta$ (by previous considerations) is less than $1/9,$ it then suffices to show that $$(k-1)\frac{80}{81}>\frac{k}{9(R-1)}.$$ As $k \ge 3$, this will be satisfied if $R>2.$ 

The mixed Ricci curvature terms can all be controlled (i.e. rendered arbitrarily small) by a suitably small choice of $c$. This means that even though the directions considered above might not be eigendirections for the Ricci tensor, the corresponding values can be viewed as estimating the genuine eigenvalues to within any given degree of accuracy. Thus the positivity of these values is sufficient for the purposes of identifying $Ric>0$ and $Ric_2>0.$ There is just one direction in which we must exercise some caution: $Ric(\partial/\partial \tau,\partial/\partial\tau)$ is only non-negative when $\tau \in [0,1/10]\cup[9/10,1].$ However the mixed Ricci curvature terms involving $\partial/\partial\tau$ all vanish, so the mixed terms do not have any impact on the Ricci curvatures of directions having a component in the $\partial/\partial\tau$ direction.

In summary then, provided $\delta<1/9$ and $R>1+(2/\delta)$ (which automatically ensures $R>2$), we will have $Ric\ge 0$ and $Ric_2>0$ if $\epsilon,c$ are chosen sufficiently small.
\bigskip

\noindent $\bullet$ {\it Around the fibre middle region.}
\medskip

We now come to the `middle' region for $\bar\psi$, corresponding to where $\theta(x)=\theta_0(x)=c\sin(\pi(x+1)/3).$ For the corresponding values of $r$ (where again the relevant range of $r$ values is dependent on $\tau$), a straightforward compuatation yields the following derivative expressions.
\begin{align*}
-\frac{1}{(1-R^{-1}\phi)^2}\frac{\bar{\psi}_{rr}}{\bar\psi}&=\frac{\pi^2}{9\Lambda^2}\cos^2\bigl(\frac{r}{R}\bigr)+\frac{\pi}{3\Lambda(R-\phi)}\sin\bigl(\frac{r}{R}\bigr)\cot\Bigl(\frac{\pi}{3}\Bigl(1+\frac{R-\phi}{\Lambda}\sin\bigl(\frac{r}{R}\bigr)\Bigr)\Bigr),\\ \\
-\frac{\bar{\psi}_{\tau\tau}}{\bar\psi}&=\frac{\pi^2}{9\Lambda^2}\phi'^2\sin^2\bigl(\frac{r}{R}\bigr)+\frac{\pi}{3\Lambda}\phi''\sin\bigl(\frac{r}{R}\bigr)\cot\Bigl(\frac{\pi}{3}\Bigl(1+\frac{R-\phi}{\Lambda}\sin\bigl(\frac{r}{R}\bigr)\Bigr)\Bigr),\\ \\
\frac{\bar{\psi}_\tau}{\bar\psi}&=-\frac{\pi}{3\Lambda}\phi'\sin\bigl(\frac{r}{R}\bigr)\cot\Bigl(\frac{\pi}{3}\Bigl(1+\frac{R-\phi}{\Lambda}\sin\bigl(\frac{r}{R}\bigr)\Bigr)\Bigr),\\ \\
\end{align*}
\begin{align*}
\frac{1-\bar{\psi}^2_\tau-(1-R^{-1}\phi)^{-2}\bar{\psi}^2_r}{\bar{\psi}^2}=&\frac{1-(\pi^2c^2/9)(\phi'^2\sin^2(r/R)+\cos^2(r/R))}{\Lambda^2c^2\sin^2\Bigl(\frac{\pi}{3}\bigl(1+\frac{R-\phi}{\Lambda}\sin(r/R)\bigr)\Bigr)}\\
&+\frac{\pi^2}{9\Lambda^2}\bigl(\phi'^2\sin^2(r/R)+\cos^2(r/R)\bigr).\\
\end{align*}

These expressions yield the Ricci curvature formulas below. In each case we analyse the formula to investigate its positivity (or otherwise).

\begin{align*}
Ric(V/\bar{\psi},V/\bar{\psi})=&(k-2)\frac{1-(\pi^2c^2/9)(\phi'^2\sin^2(r/R)+\cos^2(r/R))}{\Lambda^2c^2\sin^2\Bigl(\frac{\pi}{3}\bigl(1+\frac{R-\phi}{\Lambda}\sin(r/R)\bigr)\Bigr)}\\
&+(k-1)\frac{\pi^2}{9\Lambda^2}\bigl(\phi'^2\sin^2(r/R)+\cos^2(r/R)\bigr)\\
&+\frac{\pi}{3\Lambda}\sin\bigl(\frac{r}{R}\bigr)\cot\Bigl(\frac{\pi}{3}\Bigl(1+\frac{R-\phi}{\Lambda}\sin\bigl(\frac{r}{R}\bigr)\Bigr)\Bigr)\Bigl[\phi''+\frac{1-\phi'^2}{R-\phi}+k\phi'\frac{f'}{f}\Bigr]\\
&+\frac{\bar{\psi}^2}{f^4}(AV,AV).\\
\end{align*}
The first term in this expression will be positive provided $c<3/\pi$. The second and final terms are both non-negative.
Using the fact that $f'/f>-\delta,$ $\phi' \in [0,1]$ and $\phi''>-2$, we see that the third term is bounded below by
$$\frac{\pi}{3\Lambda}\sin(r/R)\cot\Bigl(\frac{\pi}{3}\Bigl(1+\frac{R-\phi}{\Lambda}\sin\bigl(\frac{r}{R}\bigr)\Bigr)\Bigr)\Bigl[-2-k\delta \Bigr].$$ The key point here is that for any choice of $\Lambda,$ as $c \to 0$ the first term in the expression for $Ric(V/\bar{\psi},V/\bar{\psi})$ tends to infinity, whereas the bad term above stays fixed. Therefore this Ricci curvature term is guaranteed to be positive if $c$ is sufficiently small compared to $\Lambda.$ (For the mixed Ricci curvature terms $c$ has to be sufficiently small compared to $f^{-3}$, and we will see later that in different contexts $c$ has to be small compared to $\Lambda$ for other reasons.)

Turning our attention next to the $\tau$-direction we have
\begin{align*}
Ric(\partial/\partial \tau,\partial/\partial \tau)=& \phi''\Bigl[\frac{1}{R-\phi}+(k-1)\frac{\pi}{3\Lambda}\sin\bigl(\frac{r}{R}\bigr)\cot\Bigl(\frac{\pi}{3}\Bigl(1+\frac{R-\phi}{\Lambda}\sin\bigl(\frac{r}{R}\bigr)\Bigr)\Bigr)\Bigr]\\
&+(k-1)\frac{\pi^2}{9\Lambda^2}\phi'^2\sin^2(r/R)-k\frac{f''}{f}.\\
\end{align*}
Notice that the first term in this expression is generally bad, but the other terms are good. For $\tau \in [1/10,9/10]$ we have $-kf''/f=k\delta.$ Hence by choosing $\Lambda$ (and therefore $R$) sufficiently big, we obtain positivity overall since $\delta$ is chosen independently of $\Lambda,R.$ For $\tau \in [0,1/10]\cup[9/10,1]$ we have $\phi''=0$. Therefore $Ric(\partial/\partial \tau,\partial/\partial \tau)\ge 0$ for such $\tau$. Actually, for $\tau \in [0,1/10]$ we have $\phi'>0$, so the second term in the above expression ensures the positivity of the Ricci curvature. Therefore overall we have $Ric(\partial/\partial \tau,\partial/\partial \tau)$ positive for $\tau \in[0,9/10]$, and non-negative for $\tau \in [9/10,1].$ 

In the $r$-direction we have
$$Ric\Bigl(\frac{1}{1-R^{-1}\phi}\frac{\partial}{\partial r},\frac{1}{1-R^{-1}\phi}\frac{\partial}{\partial r}\Bigr)=\frac{k\phi'}{R-\phi}\frac{f'}{f}+(k-1)\frac{\pi^2}{9\Lambda^2}\cos^2(r/R)+\frac{\phi''}{R-\phi}$$
$$+(k-1)\frac{\pi}{3\Lambda}\frac{1}{R-\phi}\sin\bigl(\frac{r}{R}\bigr)\cot\Bigl(\frac{\pi}{3}\Bigl(1+\frac{R-\phi}{\Lambda}\sin\bigl(\frac{r}{R}\bigr)\Bigr)\Bigr)[1-\phi'^2].$$
Since $1-\phi'^2\ge 0,$ we see that the final term above is non-negative. The second term is non-negative, but is not positive as $r=\pi R/2$ belongs to the middle region. Thus we are left trying to balance the bad terms $$\frac{\phi''}{R-\phi}+k\frac{\phi'}{R-\phi}\frac{f'}{f}$$ with the final non-negative term. For $\tau \in [0,1/10]$ we have $1-\phi'^2=0$ and $\phi''=0$, so we are left with a single bad term which we can bound below as follows: $$k\frac{\phi'}{R-\phi}\frac{f'}{f}>-k\frac{\delta}{R-1}.$$ Provided $\delta<1/k$ we obtain a lower bound of $-1/(R-1)$ here.
For $\tau \in [1/10,9/10]$, using $\phi''>-3/2$ and choosing $\delta<1/2k$, we can estimate the bad terms collectively by $-2/(R-1).$ As $\phi' \approx 1$ when $\tau \approx 1/10$, this is the best estimate we can give. For $\tau \in [9/10,1]$ we have Ricci positivity in the $r$-direction. Thus overall (i.e. for all $\tau \in [0,1]$), we can only say that this Ricci curvature term is bounded below by $-2/(R-1).$

Finally, we consider
\begin{align*}
Ric(X_i/f,X_i/f) =& (k-1)\frac{1-f'^2}{f^2}-\frac{f''}{f}+\frac{\phi'}{R-\phi}\frac{f'}{f}-2\frac{\bar{\psi}^2}{f^4}(A_{X_i},A_{X_i})\\
&+(k-1)\frac{\pi}{3\Lambda}\phi'\frac{f'}{f} \sin\bigl(\frac{r}{R}\bigr)\cot\Bigl(\frac{\pi}{3}\Bigl(1+\frac{R-\phi}{\Lambda}\sin\bigl(\frac{r}{R}\bigr)\Bigr)\Bigr)\\
\ge & (k-1)(1-\delta^2)-\frac{2\Lambda^2c^2}{(1-\delta)^4}(A_{X_i},A_{X_i}) \\
&-\delta\Bigl[\frac{1}{R-\phi}+(k-1)\frac{\pi}{\Lambda}\sin(r/R)\cot\Bigl(\frac{\pi}{3}\Bigl(1+\frac{R-\phi}{\Lambda}\sin\bigl(\frac{r}{R}\bigr)\Bigr)\Bigr)\Bigr], \\
\end{align*}
where we have under-estimated $f$ by $1-\delta$. For large $R,\Lambda$, the sign of $Ric(X_i/f,X_i/f)$ will agree with that of $$(k-1)(1-\delta^2)-\frac{2\Lambda^2c^2}{(1-\delta)^4}(A_{X_i},A_{X_i}).$$ Thus if $c$ is sufficiently small compared to $\Lambda$ and the $A$-tensor term, we will have $Ric(X_i/f,X_i/f)>0.$

Given that we have not been able to guarantee positive Ricci curvature, let us turn our attention to 2-positive Ricci curvature. From the above analysis, to guarantee $Ric_2>0$ for this region, we must consider the effect of adding $-2/(R-1)$ to the Ricci curvature expressions for the $X_i,$ $\tau$ and $V$ directions.

For $Ric(X_i/f,X_i/f),$ given $R,\Lambda$ sufficiently large we will not alter the positivity of this term, as for small $c$ it is still approximated by $(k-1)(1-\delta^2).$

In the case of $Ric(V/\bar{\psi},V/\bar{\psi}),$ our previous analysis of this term shows that adding the term $-2/(R-1)$ will do no damage provided $c$ is sufficiently small.

The case of $Ric(\partial/\partial \tau,\partial/\partial \tau)$ is a little more delicate. The sum of the Ricci curvature terms in the $\tau$ and $r$ directions and will be positive provided 
\begin{align*}
&-k\frac{f''}{f}+k\frac{\phi'}{R-\phi}\frac{f'}{f}+(k-1)\frac{\pi^2}{9\Lambda^2}\Bigl(\phi'^2\sin^2(r/R)+\cos^2(r/R)\Bigr) \\
&+(k-1)\frac{\pi}{3\Lambda}\frac{1}{R-\phi}\sin\bigl(\frac{r}{R}\bigr)\cot\Bigl(\frac{\pi}{3}\Bigl(1+\frac{R-\phi}{\Lambda}\sin\bigl(\frac{r}{R}\bigr)\Bigr)\Bigr)[1-\phi'^2] \\
&+\phi''\Bigl[\frac{2}{R-\phi}+(k-1)\frac{\pi}{3\Lambda}\sin\bigl(\frac{r}{R}\bigr)\cot\Bigl(\frac{\pi}{3}\Bigl(1+\frac{R-\phi}{\Lambda}\sin\bigl(\frac{r}{R}\bigr)\Bigr)\Bigr)\Bigr]>0.\\
\end{align*}
For $\tau \in [0,1/10]$ we have $\phi'=1,$ $\phi''=0$, and so the last term above is identically zero, as is the penultimate term. We are therefore left with the expression $$\frac{k}{f}\Bigl(-f''+\frac{f'}{R-\phi}\Bigr)+(k-1)\frac{\pi^2}{9\Lambda^2}.$$ The second term in the line above is positive. We claim that we can arrange for $$-f''+\frac{f'}{R-\phi} \ge 0$$ with a careful choice of $f$. Suppose that we want the downward bend of $f$ to start when $\tau=1/20$ say. For $\tau \in \bigl(\frac{1}{20},\frac{1}{20}+\frac{\delta}{100}\bigr]$, set $$f(\tau)=1-\text{exp}\Bigl(-\bigl(\tau-\frac{1}{20}\bigr)^{-2}\Bigr).$$ This smoothly extends $f(\tau)=1$ for $\tau\le 1/20$, and an easy computation shows that for $\tau \in (\frac{1}{20},\frac{1}{20}+\frac{\delta}{100}),$ $$-f''+\frac{f'}{R-\phi}=\frac{e^{-1/(\tau-(1/20))^2}}{\bigl(\tau-\frac{1}{20}\bigr)^6}\Big[4-6\bigl(\tau-\frac{1}{20}\bigr)^2-\frac{2}{R-\phi}\bigl(\tau-\frac{1}{20}\bigr)^3 \Big].$$ It is clear that for $\tau \in \bigl(\frac{1}{20},\frac{1}{20}+\frac{\delta}{100}\bigr)$, this expression is positive. For $\tau \ge \frac{1}{20}+\frac{\delta}{100},$ we can use the size of $R$ to control the term $-f''+f'/(R-\phi),$ ensuring its positivity. Thus for $\tau \in [0,1/10]$ we have $$Ric(\partial/\partial \tau,\partial/\partial \tau)+Ric\Bigl(\frac{1}{1-R^{-1}\phi}\frac{\partial}{\partial r},\frac{1}{1-R^{-1}\phi}\frac{\partial}{\partial r}\Bigr)>0.$$
For $\tau \in [1/10,9/10]$ we have $-kf''/f=k\delta$, and we can balance this term against the other terms by choosing $R,\Lambda$ sufficiently large, so as to ensure that the sum of the Ricci curvatures is positive. Finally, for $\tau \in [9/10,1]$ all terms are non-negative, and the term with factor $(1-\phi'^2)$ is strictly positive. We can therefore conclude that $$Ric(\partial/\partial \tau,\partial/\partial \tau)+Ric\Bigl(\frac{1}{1-R^{-1}\phi}\frac{\partial}{\partial r},\frac{1}{1-R^{-1}\phi}\frac{\partial}{\partial r}\Bigr)>0$$ for all $\tau \in [0,1].$ 

As noted earlier, since the mixed Ricci curvature terms are controlled by $c$, by choosing $c$ sufficiently small, we can keep the Ricci tensor eigenvalues as close as we like to the Ricci curvatures in the key directions above, and hence the mixed curvature terms will not affect $Ric_2>0$.
\bigskip

\noindent $\bullet$ {\it Smoothing issues.}
\medskip

It remains to discuss what happens to the Ricci curvature over the regions of smoothing, when $\theta_0$ is smoothed to the function $\theta$. It suffices to consider the smoothing of $\theta_0(x)$ at the point $x=x_1,$ as the symmetry of $\theta_0$ about $x=1/2$ yields the same outcome at the other smoothing point.

The smoothing region can be kept arbitrarily small, and hence in this region we can assume $\theta \approx \theta_0,$ with $\theta'$ decreasing rapidly between the values of $\theta_0'$ on either side of the region. As a consequence, $-\theta''$ experiences a large positive `spike'. We need to consider the effect this has on the terms of the Ricci curvature formulas depending on $\theta$ (i.e. the terms involving $\bar\psi$). The relevant expressions are treated individually below. Note that smoothing $\theta_0$ to $\theta$ around $x=x_1$ is equivalent to smoothing $\psi$, viewed as a function on the $(r,\tau)$-plane, in a neighbourhood of the curve determined by $$\frac{R-\phi(\tau)}{\Lambda}\sin(r/R)=x_1.$$ 

In the case of $\bar{\psi}_r,$ we have a linear dependence on $\theta'$, so $\bar{\psi}^2_r$ decreases over the smoothing.

In the case of $\bar{\psi}_\tau,$ we also have a linear dependence on $\theta'$, so $\bar{\psi}_\tau/\bar{\psi}$ and $\bar{\psi}^2_\tau$ both decrease in magnitude over the smoothing. Note that these terms are both bad for the curvature, so a reduction in the size of these terms is therefore good.

In the case of $-\bar{\psi}_{\tau\tau}/\bar\psi,$ we have a good term (depending linearly on $\theta''$, assuming the variation in $\theta$ across the smoothing region is negligible), and a bad term (which likewise depends essentially linearly on $\theta'$). Thus the bad term gets smaller, and the good term spikes as it transitions between its values on either side of the smoothing. So overall, this quantity is better over the smoothing region than at either side.

In the case of $-(1-R^{-1}\phi)^{-2}\bar{\psi}_{rr}/\bar\psi,$ we have two good terms: one effectively linear in $\theta''$ and the other effectively linear in $\theta'$. The first of these improves over the smoothing, but the second diminishes. However this does not impact on $Ric_2>0$ because of the following observation.

If we have $Ric_2>0$ on either side of the smoothings, and the various terms which make up the curvature formulas either stay approximately the same or vary approximately linearly over the smoothing intervals, then the $Ric_2>0$ condition 
will be preserved provided we can control the deviation from being constant or varying linearly. In our case we can control these factors by limiting the size of the smoothing region. We also have a term ($\theta''$) which is behaves in a way which is better than if it varied linearly.

Indeed, checking all the Ricci curvature expressions (including the mixed terms), we see that all these curvatures either improve, or in the case of the $r$-direction certainly get no worse, over the smoothing regions. Hence $Ric_2>0$ is preserved if these regions are sufficiently small.
\medskip

We summarize the all the above considerations in the following
\begin{proposition}\label{second_c_restriction}
Given $\delta$ with $0<\delta<\min\{1/2k,1/9\}$, there exists $R_0=R_0(\delta)$ with the following property. For any $R>R_0,$ there exists $c_0=c_0(\delta,R,\nabla),$ such that after setting $\Lambda=5R/2$, any choice of $c \in (0,c_0)$ results in the metric $g_\phi$ having $Ric_2>0.$
\end{proposition}
\bigskip

\centerline{\it Completing the surgery on $W$}
\bigskip\medskip

To complete the surgery on $W$, it remains to `cap off' $W$ with a disc $D^{k+1}$ (having radial parameter $\tau$), and extend the metric $g_\phi$ to a $Ric_2>0$ metric over $D^{k+1} \times S^k.$

Let us focus for the moment on the normal sphere metric induced by $g_\phi$ corresponding to $\tau \approx 1.$ This is, for $r \in [0,R\pi],$ 
\begin{equation}
\label{induced_metric}
\Bigl(1-\frac{1}{2R}\Bigr)^2dr^2+\Lambda^2\theta^2\Bigl(\frac{R-\frac{1}{2}}{\Lambda}\sin\bigl(\frac{r}{R}\bigr)\Bigr)ds^2_{k-1}. 
\end{equation}

Recall that this sphere metric is independent of $\tau$ when $\tau \approx 1.$

We now want to think of this metric as a scaling of a `model' fibre metric $\omega$, which is a single warped product metric over a fixed interval $[0,\pi]$. We set
\begin{equation}
\label{omega}
\omega=d\eta^2+\frac{\Lambda^2}{(R-\frac{1}{2})^2}\theta^2\Bigl(\frac{R-\frac{1}{2}}{\Lambda}\sin(\eta)\Bigr)ds^2_{k-1}
\end{equation}
for $\eta \in [0,\pi].$ The metric (\ref{induced_metric}) is then $$\Bigl(R-\frac{1}{2}\Bigr)^2\omega,$$ where $r/R$ has been replaced by $\eta$. For convenience we will set $$\bar{\theta}(\eta):=\Bigl(\frac{\Lambda}{R-\frac{1}{2}}\Bigr)\theta\Bigl(\frac{R-\frac{1}{2}}{\Lambda}\sin(\eta)\Bigr),$$ so $$\omega=d\eta^2+\bar{\theta}^2(\eta)ds^2_{k-1}.$$

In order to estimate the Ricci curvatures of $\omega$ we make the derivative computations below.
$$\bar{\theta}'(\eta)=\cos(\eta)\theta'\Bigl(\frac{R-\frac{1}{2}}{\Lambda}\sin(\eta)\Bigr);$$
$$\bar{\theta}''(\eta)=\Bigl(\frac{R-\frac{1}{2}}{\Lambda}\Bigr)\cos^2(\eta)\theta''\Bigl(\frac{R-\frac{1}{2}}{\Lambda}\sin(\eta)\Bigr)-\sin(\eta)\theta'\Bigl(\frac{R-\frac{1}{2}}{\Lambda}\sin(\eta)\Bigr).$$
Therefore
$$-\frac{\bar{\theta}''}{\bar\theta}=-\Bigl(\frac{R-\frac{1}{2}}{\Lambda}\Bigr)^2\cos^2(\eta)\frac{\theta''}{\theta}\Bigl(\frac{R-\frac{1}{2}}{\Lambda}\sin(\eta)\Bigr)+\Bigl(\frac{R-\frac{1}{2}}{\Lambda}\Bigr)\sin(\eta)\frac{\theta'}{\theta}\Bigl(\frac{R-\frac{1}{2}}{\Lambda}\sin(\eta)\Bigr).$$
Observe that a lower bound for this last expression gives a lower bound for the Ricci curvature of $\omega$ (as $|\theta'|\le 1$). 

Notice that the first term in the expression for $-\bar{\theta}''/\bar{\theta}$ is positive for $\eta \neq \pi/2,$ and the second term is positive for $\eta \neq 0,\pi$. As $\Lambda=5R/2$ we see that provided $R\ge 3$ (which we can assume), $$\frac{R-\frac{1}{2}}{\Lambda}\in\Bigl[\frac{1}{3},\frac{2}{5}\Bigr).$$ Thus $$\theta'\Bigl(\frac{R-\frac{1}{2}}{\Lambda}\sin(\eta)\Bigr)>0 \,\, \text{ for all } \eta \in [0,\pi].$$

In order to obtain a lower bound for the whole expression,
we split the $\eta$-interval $[0,\pi]$ into $[0,\pi/3]\cup [2\pi/3,\pi]$ and $[\pi/3,2\pi/3].$

On the union of intervals above, we have $\cos^2(\eta)>\cos^2(\pi/3)=1/4.$ Thus we obtain a lower bound here for $-\bar{\theta}''/\bar\theta$ of $$-\frac{\bar{\theta}''}{\bar\theta} >\frac{1}{9}\times\frac{1}{4}\times\Bigl(-\frac{\theta''}{\theta}\Bigr).$$ 
By Lemma \ref{theta_est} we have $-\theta''/\theta>1$ and so we obtain a lower bound $$-\frac{\bar{\theta}''}{\bar\theta} >\frac{1}{36}$$ on the union of intervals.

On $[\pi/3,2\pi/3]$ we consider the second term in the expression for $-\bar{\theta}''/\bar\theta$. Note first that by our choice $\Lambda,R$, the argument of $\theta,$ namely $\frac{R-\frac{1}{2}}{\Lambda}\sin(\eta)$, always lies between $\frac{\sqrt 3}{6}$ and $2/5$ say, over this interval. Thus the argument lies well inside the interval between the smoothing regions for $\theta_0$. In this middle region we have a lower bound for this second term of $$\frac{1}{3}\sin\bigl(\frac{\pi}{3}\bigr)\frac{\theta'}{\theta}\Bigl(\frac{R-\frac{1}{2}}{\Lambda}\sin(\eta)\Bigr)=\frac{1}{2\sqrt 3}\frac{\theta'}{\theta}\Bigl(\frac{R-\frac{1}{2}}{\Lambda}\sin(\eta)\Bigr),$$ and we also know that $$\frac{\theta'}{\theta}(x)=\frac{\pi}{3}\cot\bigl(\frac{\pi}{3}(x+1)\bigr).$$ Since the cotangent function is decreasing we have
\begin{align*}
\frac{\theta'}{\theta}\Bigl(\frac{R-\frac{1}{2}}{\Lambda}\sin(\eta)\Bigr)&>\frac{\pi}{3}\cot\bigl(\frac{\pi}{3}\times \frac{7}{5}\bigr) \\
&> 0.11.\\
\end{align*}
Therefore over this middle interval, the second term in the expression for $-\bar{\theta}''/\bar\theta$ is bounded below by $$\frac{1}{2\sqrt 3} \times 0.11 \approx 0.0318>\frac{1}{36}.$$ So overall (i.e. taking the whole of $[0,\pi]$ into account), we obtain the lower bound $$-\frac{\bar{\theta}''}{\bar\theta}>\frac{1}{36}.$$

Note in particular that we have a lower bound for $-\bar{\theta}''/\bar\theta$ which is independent of all other constants introduced so far (namely $c,\Lambda,R,\delta$). (For the record, the relationship between $\Lambda$ and $R$ was chosen so as to give such a lower bound.) Using the single warped product formulas we immediately obtain
\begin{proposition}
\label{1/36}
For all unit vectors, the Ricci curvatures of $\omega$ are bounded below by $1/36$ (independent of all choices used in the construction of $\omega$).
\end{proposition}

We now face a situation similar to that encountered in the `transition phase' of the metric construction on $S^{2k+1}$ in Section \ref{sec:metrics}: the metric $g_\phi$ has the wrong symmetries to extend over $D^{k+1} \times S^k$ without altering the metric $\omega$. Specifically, just as in Section \ref{sec:metrics}, we need to deform $\omega$ to a round metric, and then use the increased symmetry to flatten the connection, so we end up with a product metric we can easily cap off. Even though $\omega$ is not isometric to $ds^2+\psi^2(s)ds^2_{k-1}$ which was rounded in Section \ref{sec:metrics} (see Lemma \ref{can_var} onwards), essentially the same arguments work here. However for our current purposes we will require a more detailed analysis.

Note that the submersion metrics below are all relative to the projection $\pi:S^k_a \times S^k_b \to S^k_a$.

Given some $\lambda>0,$ we will join the metric $\nu(\lambda^2 ds^2_k,\rho^2\omega,\nabla)$ to $\nu(\lambda^2 ds^2_k,\rho^2 ds^2_k,\nabla_{triv})$ along a path of Ricci positive metrics, for some choice of $\rho>0.$

Clearly, any upper bound for $\rho$ will depend on $\nabla$ and $\lambda$. We claim that this is the only thing on which $\rho$ depends, in that the upper bound for $\rho$ is independent of the choices (of $\Lambda,$ $R=2\Lambda/5,$ and $c$) used to construct $\omega$.

For $v \in [0,1],$ set $$H(\eta,v):=(1-\beta(v))\bar{\theta}(\eta)+\beta(v)\sin(\eta),$$ where $\beta$ is the bump function introduced before Lemma \ref{H}. (Note that $H$ here is unrelated to the function $H(s,x)$ which appeared in Section \ref{sec:metrics}.)
Consider the path of metrics $P(v)$ given by $$P(v)=d\eta^2+H^2(\eta,v)ds_k^2.$$ Observe that since $H$ is a convex combination of odd functions and the $\eta$-derivatives $H_\eta(0,v)=H_\eta(\pi,v)=1,$ $P(v)$ is a smooth metric on $S^k$ for all $v \in [0,1].$ Moreover $P(v)$ is a smooth path linking $\omega$ to $ds^2_k.$

\begin{proposition}\label{bar_rho_1}
Given $\lambda>0$, with the path of metrics $P(v),$ $v \in [0,1],$ as above, there exists $\kappa_1=\kappa_1(\nabla,\lambda) \in (0,1)$ such that for all $\rho \in (0,\kappa_1),$ the submersion metric $$\nu\Bigl(\lambda^2 ds^2_k,\rho^2P(v),\nabla\Bigr)$$ has positive Ricci curvature for all $v \in [0,1],$ irrespective of the choices involved in the construction of $\omega$.
\end{proposition}

\begin{proof} 
The Ricci positivity of the submersion metric in question depends on the base metric, the connection, and the fibre metric. More specifically, it depends on the Ricci curvature of the base, the $A$-tensor arising from the connection, and the lower bound for the Ricci curvatures of the fibre. (See \cite[\S9G]{Be}.) The result follows easily from Lemma \ref{can_var} if we can establish that $Ric(P(v))$ has a lower bound independent of the choices made to construct $\omega.$ 
From the single warped product formulas we see that for each $v$, a lower bound for $Ric(P(v))$ is given by $-H_{\eta\eta}/H$. It therefore suffices to show that $-H_{\eta\eta}/H$ is positive and independent of the choices made to define $\omega$.

Observe that for $v=1,$ $$-\frac{H_{\eta\eta}}{H}=-\frac{\sin''(\eta)}{\sin(\eta)}=1,$$ and recall that for $v=0,$ $$-\frac{H_{\eta\eta}}{H}=-\frac{\bar{\theta}''}{\bar\theta}>\frac{1}{36}$$ (from Proposition \ref{1/36} above). We must therefore check that for $v \in (0,1)$, $-H_{\eta\eta}/H$ maintains a positive lower bound independent of all choices. Below we consider the cases $\eta=0,\pi$ and $\eta \in (0,\pi)$ separately.

Suppose that $\eta=0$. (The case $\eta=\pi$ is essentially the same by symmetry.) As $H_{\eta\eta}(0,v)=H(0,v)=0,$ we investigate $-H_{\eta\eta}/H$ at $(0,v)$ using l'H\^opital's rule. This gives $$-\frac{H_{\eta\eta}}{H}(0,v)=-\lim_{\eta \to 0}\frac{H_{\eta\eta\eta}}{H_\eta}(0,v).$$ Now it is easily checked that $H_\eta(0,v)=1,$ and therefore $$-\frac{H_{\eta\eta}}{H}(0,v)=-\lim_{\eta \to 0}H_{\eta\eta\eta}(0,v).$$ But $H_{\eta\eta\eta}$ is a convex combination of $\bar{\theta}'''(\eta)$ and $\sin'''(\eta)$, and by l'H\^opital again we also have $$-\bar{\theta}'''(0)=-\lim_{\eta \to 0}\frac{\bar{\theta}''}{\bar\theta}(\eta)>\frac{1}{36}.$$ Therefore $$-\frac{H_{\eta\eta}}{H}(0,v)> (1-\beta)\frac{1}{36}+\beta.1 \ge \frac{1}{36}.$$ For $\eta \neq 0,\pi$, the result follows immediately from the next lemma.
\end{proof}

\begin{lemma}
Given non-zero real-valued functions $A,B$ defined on the same interval, with $-A''/A$ and $-B''/B$ strictly positive, we have at any point in the domain interval $$\frac{(1-\beta)(-A'')+\beta(-B'')}{(1-\beta)A+\beta B}\ge \min\Bigl\{-\frac{A''}{A},-\frac{B''}{B}\Bigr\}$$ for any $\beta\in [0,1].$
\end{lemma}

\begin{proof}
Suppose without loss of generality that at some point in the domain interval we have $$-\frac{A''}{A} \ge -\frac{B''}{B},$$ or equivalently $$-A''B \ge -B''A.$$ Then
\begin{align*}
\frac{(1-\beta)(-A'')+\beta(-B'')}{(1-\beta)A+\beta B}&=\frac{(1-\beta)(-A'')+\beta(-B'')}{(1-\beta)A+\beta B}\frac{AB}{AB}\\ \\
&\ge\frac{(1-\beta)(-B''A)A+\beta(-B''A)B}{(1-\beta)A(AB)+\beta B(AB)}\\ \\
&=-\frac{B''A}{AB}\frac{(1-\beta)A+\beta B}{(1-\beta)A+\beta B}\\ \\
&=-\frac{B''}{B}.\\
\end{align*}
Hence the result.
\end{proof}

Given $\lambda>0,$ it will be convenient to choose $\kappa_1$ to be supremum of the permissable values in $(0,1)$ guaranteed by Proposition \ref{bar_rho_1}. 

Having deformed the fibre metric to be round, we next untwist the connection. Let $\nabla(v),$ $v \in [1,2]$ say, be any path of principal connections joining $\nabla$ to $\nabla_{triv}$. It follows immediately from Lemma \ref{can_var} that there exists $\kappa_2=\kappa_2(\nabla(v),\lambda) \in (0,1)$ such that $$\nu(\lambda^2 ds^2_k, \rho^2 ds^2_k,\nabla(v))$$ has positive Ricci curvature for each $v \in [1,2]$ and any $\rho \in (0,\kappa_2).$ Again, we will assume that the chosen value of $\kappa_2$ is the supremum of the possible values within $(0,1).$ 

Set $\bar{\rho}=\min\{\kappa_1,\kappa_2\}.$ Thus $\bar{\rho}$ depends continuously on $\lambda$. Concatenating the paths above and using Lemma \ref{transition} gives

\begin{lemma}
\label{tt}
Given $\lambda>0$ and $\bar\rho>0$ as above, for all $\rho\in(0,\bar\rho)$ there is a $Ric_2>0$ metric on $I \times S^k \times S^k$ (with the interval $I$ parametrized by $u$), such that close to one boundary component the metric agrees with $du^2+\nu(\lambda^2 ds^2_k,\rho^2\omega,\nabla),$ and close to the other it agrees with $du^2+\lambda^2ds^2_k+\rho^2ds^2_k.$
\end{lemma}

\begin{corollary}\label{cap_off}
With $\lambda,\bar{\rho}$ as before, for any $\rho \in (0,\bar\rho)$ there is a $Ric_2>0$ metric on $D^{k+1}\times S^k$ such that close to the boundary the metric agrees with $du^2+\nu(\lambda^2 ds^2_k,\rho^2\omega,\nabla),$ and in a neighbourhood of the centre takes the form $$du^2+\sin^2(u)ds^2_k+\rho^2ds^2_k,$$ where we are now taking $u$ to be the radial parameter in $D^{k+1},$ with the centre of $D^{k+1}$ corresponding to $u=0.$
\end{corollary}

\begin{proof}
We need a double warped product metric of the form $$du^2+F^2(u)+\rho^2ds^2_k$$ such that for $u \approx 0$ we have $F(u)=\sin(u),$ and for $u \ge u_0,$ some $u_0>0$ we have $F(u)=\lambda.$ From the double warped product Ricci curvature formulas it is clear that any concave down function $F(u)$ satisfying the above boundary conditions will do.
\end{proof}

It is clear that any $\lambda>0$ will work here. However looking at the wider picture, under suitable concatenation of parameters, we will need the function $F$ in the proof of the Corollary above to smoothly join the function $f$ discussed previously. This will certainly require $\lambda<1$. Without loss of generality, let us choose $\lambda=1/2.$ Consider the corresponding value $\bar{\rho}$, and fix $\rho \in (0,\bar{\rho}).$
\medskip

To complete the $Ric_2>0$ surgery it remains to show how to transition between $g_\phi$ and $d\tau^2+\nu(\lambda^2ds^2_k,\rho^2\omega,\nabla)=d\tau^2+\nu(\frac{1}{4}ds^2_k,\rho^2\omega,\nabla).$

Near the relevant boundary, $g_\phi$ takes the form $d\tau^2+\nu(f^2ds^2_k,(R-\frac{1}{2})^2\omega,\nabla).$ Thus to make the transition we must do two things: change the base scaling from $f\ (\approx 1)$ to $\lambda=1/2$, and change the fibre scaling from $(R-\frac{1}{2})$ to $\rho.$ It will be convenient to work in the opposite direction to that of increasing $\tau$: to this end consider a parameter $u\ge u_0,$ some $u_0>0$, extending the $u$-parameter in Corollary \ref{cap_off}, with the boundary of $D^{k+1}$ in the Corollary corresponding to $u=u_0$. Moreover, suppose that for some $u_0+N>u_0$ we have $u-u_0-N=1-\tau$ for $u \in [u_0+N,u_0+N+1]$. We want the function $F(u)$ to `extend' $f(\tau)$ in the sense that $F(u)=f(1-(u-u_0-N))$ for $u \in [u_0+N,u_0+N+1]$.

We will therefore consider metrics of the form $$du^2+\nu(F^2(u)ds^2_k,h^2(u)\omega,\nabla)$$ for $u \ge u_0,$ and we will demand that $h(u)=\rho$ for $u \in [u_0,u_0+(1/10)]$, and $F(u)=\lambda=1/2$ for $u \in [u_0,u_0+1].$ 

A straightforward calculation (see the Appendix) yields the following Ricci curvature formulas for this metric:
\begin{lemma}\label{end_ric} 
The Ricci curvatures of the metric $du^2+\nu(F^2(u)ds^2_k,h^2(u)\omega,\nabla)$ are given by
\begin{align*}
Ric(V/h\bar\theta,V/h\bar\theta)=&\frac{1}{h^2}\Bigl(-\frac{\bar{\theta}''}{\bar{\theta}}+(k-2)\frac{1-\bar{\theta}'^2}{\bar{\theta}^2}\Bigr)-(k-1)\frac{h'^2}{h^2}-\frac{h''}{h}\\
& -k\frac{h'F'}{hF}+\frac{h^2\bar{\theta}^2}{F^4}(AV,AV);\\
Ric(X_i/F,X_i/F)=&-\frac{F''}{F}+(k-1)\frac{1-F'^2}{F^2}--k\frac{h'F'}{hF}-\frac{2\bar{\theta}^2h^2}{F^4}(A_{X_i},A_{X_i});\\
Ric(\partial/\partial u,\partial/\partial u)=&-k\frac{F''}{F}-k\frac{h''}{h};\\
Ric(\frac{1}{h}\frac{\partial}{\partial \eta},\frac{1}{h}\frac{\partial}{\partial \eta})=&-(k-1)\frac{\bar{\theta}''}{\bar\theta}\frac{1}{h^2}-(k-1)\frac{h'^2}{h^2}-\frac{h''}{h}-k\frac{h'F'}{hF};\\
Ric(X_i/F,V/h\bar\theta)=&-\frac{h\bar\theta}{F^3}(\check{\delta}A(X_i),V).\\
\end{align*}
Here, as before $V$ is a unit vector for $(S^{k-1},ds^2_{k-1})$ with $S^{k-1}$ a cross-section of the fibre, and the $A$-tensor terms are those for the `comparison' sphere bundle $(\mathcal{B},g_B)$. All other mixed Ricci curvature terms vanish.
\end{lemma}

Some estimation of terms in the above formulas will be useful later. 
\begin{lemma} \label{lower_bounds}
The following lower bounds hold for all $\eta \in [0,\pi]$:
\begin{enumerate}
\item $-(k-1)\frac{\bar{\theta}''}{\bar{\theta}}>(k-1)/36$;
\item $-\frac{\bar{\theta}''}{\bar{\theta}}+(k-2)\frac{1-\bar{\theta}'^2}{\bar{\theta}^2}>(k-1)/36.$
\end{enumerate}
\end{lemma}

\begin{proof}
The first inequality above was established in the proof of Proposition \ref{1/36}. The second inequality will follow if we can show that $(1-\bar{\theta}'^2)/\bar{\theta}^2>1/36.$ Directly from the definition of $\bar{\theta}$ we see that 
\begin{align*}
\frac{1-\bar{\theta}'^2}{\bar{\theta}^2}&=\Bigl(\frac{R-\frac{1}{2}}{\Lambda}\Bigr)^2\Bigg[\frac{1-\cos^2(\eta)\theta'^2\Bigl(\frac{R-\frac{1}{2}}{\Lambda}\sin(\eta)\Bigr)}{\theta^2\Bigl(\frac{R-\frac{1}{2}}{\Lambda}\sin(\eta)\Bigr)}\Bigg] \\
&\ge \Bigl(\frac{R-\frac{1}{2}}{\Lambda}\Bigr)^2\Bigg[\frac{1-\theta'^2\Bigl(\frac{R-\frac{1}{2}}{\Lambda}\sin(\eta)\Bigr)}{\theta^2\Bigl(\frac{R-\frac{1}{2}}{\Lambda}\sin(\eta)\Bigr)}\Bigg] \\
& \ge \Bigl(\frac{R-\frac{1}{2}}{\Lambda}\Bigr)^2, \\
\end{align*}
where the last line is a consequence of Lemma \ref{theta_est}. But $(R-\frac{1}{2})/\Lambda=(2R-1)/5R,$ which (assuming $R \ge 3$ ) lies between $1/3$ and $2/5.$ As $(1/3)^2>1/36,$ the claim follows.
\end{proof}

Our next task is to describe the construction of the functions $h$ and $F$. Let $\beta:{\mathbb R} \to {\mathbb R}$ now denote the obvious extension of the bump function introduced before Lemma \ref{H}, i.e. $\beta(x)=0$ for $x \le 0,$ $\beta(x)=1$ for $x \ge 1$, $\beta'(x) \ge 0$ for all $x \in {\mathbb R}.$ For constants $s_1,s_2>0$ (to be chosen later), set $$h(u)=\rho+s_1\int_{u_0}^u \beta(x-u_0-(1/10))\,dx,$$ $$F(u)=\frac{1}{2}+s_2\int_{u_0}^u \beta(x-u_0-1)\, dx.$$ Notice that these functions satsify the requirements detailed above, i.e. $h(u)=\rho$ for $u \in [u_0,u_0+(1/10)]$ and $F(u)=1/2$ for $u \in [u_0,u_0+1].$ Moreover, for $u \ge u_0+(11/10)$ respectively $u \ge u_0+2$, $h$ and $F$ have constant slopes $s_1$ respectively $s_2$.

Given the constant slope behaviour observed above, the following result will be crucial:
\begin{lemma}\label{const_slope}
For any $a,b,C,C'>0,$ the inequalities $$\frac{C}{(s_1x+a)^2}>\frac{s_1s_2}{(s_1x+a)(s_2x+b)},$$ $$\frac{C'}{(s_2x+b)^2}>\frac{s_1s_2}{(s_1x+a)(s_2x+b)}$$ hold simultaneously for all $x \ge 0$, provided $$0<s_1<\min\{1,\sqrt C, bC/a\}$$ and $$0<s_2<\min\{1,\sqrt{C'},aC'/b\}.$$
\end{lemma}

\begin{proof}
A straightforward rearrangement of the desired inequalities shows that we need $$s_1<\min\{\sqrt C, bC/(s_2a)\} \text{ and }s_2<\min\{\sqrt{C'},aC'/(s_1b)\}.$$ Restricting $s_1,s_2$ to the interval $(0,1)$ then means that we can replace $s_1$ and $s_2$ in the above minima by 1.
\end{proof}

\begin{corollary}
\label{ttt}
There exist $s_1=s_1(\rho,k)$ and $s_2=s_2(\rho,k,s_1)$ in the interval $(0,1),$ with $s_2<s_1$, such that for the resulting functions $F(u),h(u)$ defined above, the metric $$du^2+\nu(F^2(u)ds^2_k,h^2(u)\omega,\nabla)$$ has $Ric_2>0$ (and $Ric \ge 0$) for $u \in [u_0,u_0+N]$, for any given $N>2$, provided $c$ is chosen sufficiently small depending on the other parameters.
\end{corollary}

\begin{proof}
We analyse the formulas in Lemma \ref{end_ric}, given the form of the functions $h$ and $F$. 

We begin by observing that $$\bar{\theta}_{max}=\frac{\Lambda}{R-\frac{1}{2}}c.$$ Since $\Lambda=5R/2$ we see that assuming $R \ge 3,$ the maximum value of $\Lambda/(R-\frac{1}{2})$ is 3. We can therefore estimate $\bar{\theta}_{max}$ by $3c,$ and so by choosing $c$ sufficiently small, we can render the $A$-tensor terms as small as we like over any given compact interval. Since we are working on the compact interval $u \in [u_0,u_0+N],$ this means in effect that we can ignore the $A$-tensor terms appearing in the curvature formulas: if setting the $A$-tensor terms temporarily to zero results in $Ric_2>0$, then by the openness of this condition, we can simply choose $c$ sufficiently small so that the curvature condition is preserved when the $A$-tensor terms are re-introduced. We will therefore {\it assume that the $A$-tensor expressions all vanish}, and consider any $u \ge u_0.$ 
  

It follows trivially from Lemma \ref{end_ric} that for $u \in [u_0,u_0+(1/10)]$ we have $Ric_2>0$ and $Ric \ge 0$. For $u \in [u_0+(1/10),u_0+(11/10)],$ the function $h$ bends upwards to achieve a constant slope $s_1$. By the openness of the $Ric_2>0$ condition, it is clear that for $s_1$ sufficiently small, the metric will continue to have $Ric_2>0$ over this interval. Temporarily fix such a value of $s_1$. 

Now repeat the same argument for $F$ over the interval $[u_0+1,u_0+2]$: if $s_2$ is small enough, by openness we will have $Ric_2>0$ over this interval. Temporarily fix such an $s_2$, subject to the condition that $s_2<s_1$.

For $u \ge u_0+2$, both $h$ and $F$ have constant slope. To show that the $Ric_2>0$ condition continues to hold for all such $u$, we use Lemma \ref{const_slope} in conjunction with Lemma \ref{lower_bounds}. By the latter Lemma, it suffices to ensure the positivity of the following expressions, arising from the $V$ and $\eta$ directions in the first case, and the $X_i$ direction in the second:
$$(k-1)\frac{1}{h^2}\Bigl(\frac{1}{36}-s_1^2\Bigr)-k\frac{s_1s_2}{hF};$$
$$(k-1)\frac{1-s_2^2}{F^2}-k\frac{s_1s_2}{hF}.$$
If we assume that $s_1 \in (0,1/\sqrt{72})$ and $s_2 \in (0,1/\sqrt 2),$ it then suffices to show that $$\frac{k-1}{72h^2}-k\frac{s_1s_2}{hF}>0;$$
$$\frac{k-1}{2F^2}-k\frac{s_1s_2}{hF}>0.$$
Setting $x=u-u_0-2,$ $a=h(u_0+2),$ $b=F(u_0+2),$ $C=\frac{k-1}{72k}$ and $C'=\frac{k-1}{2k}$, Lemma \ref{const_slope} gives conditions on $s_1,s_2$ which guarantee this desired positivity. Given that $\int_0^1 \beta(x)\,dx<1,$ we see that $\rho<a<\rho+1$ and $1/2<b<3/2.$ Thus in order to comply with the restrictions from Lemma \ref{const_slope} it suffices to ensure that $$0<s_1<\min\Bigl\{1,\sqrt C,\frac{C}{2(\rho+1)}\Bigr\}$$ and $$0<s_2<\min\Bigl\{1,\sqrt{C'},\frac{2\rho C'}{3}\Bigr\}.$$
If our choices of $s_1,s_2$ already satisfy these conditions, then we are done. If $s_1$ satisfies the condition but $s_2$ does not, simply choose a smaller value of $s_2$: re-choosing $s_2$ in this way has no other consequences. If, on the other hand, $s_1$ does not satisfy the above inequality, we must first re-choose it so as to comply with this restriction. Even if $s_2$ satisfies its inequality, it might happen that the condition $s_2<s_1$ is violated by the smaller choice of $s_1$. Either way, we might be forced us to choose a lower value for $s_2$, but again, we are free to do this if necessary. Thus in all circumstances, we can find values of $s_1,s_2$ as in the statement of the Corollary.
\end{proof}

We can now smoothly glue the metric $du^2+\nu(F^2(u)ds^2_k,h^2(u)\omega,\nabla)$ defined on $[u_0,u_0+N] \times S^k \times S^k$ to the metric on $D^{k+1}(u_0) \times S^k$ from Corollary \ref{cap_off}, to give a metric on $D^{k+1}(u_0+N) \times S^k.$ Our final surgery task is to smoothly join this to $g_\phi$ within $Ric_2>0.$

Before explaining how to complete the surgery, we are now in a position to to choose all the constants and functions needed for our constructions.

We begin by choosing a suitable value for $\rho_0$ in Proposition \ref{V_surgery}: this has no other consequences. We have already chosen $\rho$ above, and by Corollary \ref{ttt} we can choose and fix a suitable value for $s_1$, and a corresponding (but temporary) value for $s_2$. Next, pick a value for $\delta$ such that 
\begin{equation}
\label{delta}
\delta \in \bigl(0,\min\{s_2,1/(2k),1/9\}\bigr).
\end{equation}
Having selected this $\delta$, choose and fix a suitable function $f(\tau)$ (see the criterion after the statement of Proposition \ref{ric(g_phi)}) for $\tau \in [0,1]$. Notice that $f'(1)>-\delta.$ We will assume without loss of generality that $f(1)>1/2.$

Now re-choose and fix $s_2$ so that $s_2=-f'(1).$ This is a smaller value then was temporarily selected above, but notice that choosing $s_2$ smaller has no consequences for Corollary \ref{ttt}. The key point here is that setting $s_2=-f'(1)$ facilitates the creation of a $C^1$-join where $f$ and $F$ meet.

With $s_1$ and $s_2$ selected, we have now in effect determined the functions $h(u)$ and $F(u).$

Next fix $R>R_0(\delta)$ as in Proposition \ref{second_c_restriction}. (As $R>1+(2/\delta)$ and $\delta<1/9$, we see that $R>19.$) Then set $\Lambda=5R/2.$

Now choose $T>6R$ in accordance with Corollary \ref{T_and_rho_1}, and then select $\rho_1$ as indicated in the same result.  We can then fix functions $p(t),q(t)$ as introduced at the start of Section \ref{sec:metrics}, as both $\rho_0$ and $\rho_1$ have been chosen.

As $h'>F'$ (i.e. $s_1>s_2$) for $u \ge u_0+2,$ $\lim_{u \to \infty}h(u)/F(u)=\infty.$ In particular there exists $N>2$ such that $$\frac{R-\frac{1}{2}}{f(1)}=\frac{h}{F}(u_0+N).$$ (We have $R-\frac{1}{2}>37/2$ by the above, and $f(1)<1.$ Hence the ratio on the left-hand side above is $>13/2.$ In contrast, $h/F$ takes the value $2\rho<2$ when $u=u_0.$ Since $h'<1,$ it follows that $N>2$.) Fix this value of $N$.

Finally, we choose $c$ sufficiently small so as to simultaneously satisfy Proposition \ref{first_c_restriction}, Proposition \ref{second_c_restriction} and Corollary \ref{ttt}. This completes the selection of data needed for our construction.
\smallskip

In order to complete the surgery we globally rescale the metric on $D^{k+1}(u_0+N) \times S^k.$ The rescaling factor is chosen so that the function $\bar h$ obtained by rescaling $h$ has the value $R-\frac{1}{2}$ at the boundary, which (by the choice of $N$) forces the value of $\bar{F},$ the rescale of $F$, to have value $f(1).$

We make a small concave down adjustment to $\bar h$ near the boundary to render it constant with value $R-\frac{1}{2}$, making an arbitrarily small $C^0$ change to the function over the deformation. As is clear from the Ricci curvature formulas above, this only has a positive effect on the Ricci curvature.

Notice that the global rescale does not affect the scaling functions' first derivative behaviour at the boundary. In particular, we obtain a $C^1$-join between $F$ and $f$ at $u=u_0+N,$ equivalently $\tau=1.$ Thus we can glue the rescaled metric, which of course still has $Ric_2>0$, to $g_\phi$ to create a $C^1$-join. Standard arguments now apply to make the metric $C^2$ by smoothing the $C^1$-join between $f$ and $\bar{F}.$ This can be done keeping the values of the (common) function and its first derivatives approximately the same, with second derivatives interpolating approximately linearly between their original values either side of the join. As the Ricci curvature is linear in the second derivatives of the metric, we can clearly smooth preserving the $Ric_2>0$ condition.

This completes the metric construction for the surgery on $W$, and establishes that this operation can be performed within $Ric_2>0$ as claimed.
\medskip

In order to prove Theorem \ref{thm:main} we still have to consider the possibility of surgeries performed with alternative trivialisations of the normal bundle. This applies to both the surgery on $W$ as considered above, and also the surgery on $Z$ addressed in Proposition \ref{V_surgery}.

The issue of different trivialisations is discussed in detail in \cite[\S1]{Wr2}, but see also the synopsis in \cite[\S3]{Wr4}. For the benefit of the reader, we will outline below the basic arguments needed.

Consider again the situation of Corollary \ref{cap_off}, which shows how the original surgery on $W$ can be metrically capped-off. Specifically, we have a metric on the `cap' $D^{k+1} \times S^k$ which near its centre is isometric to the product of a unit round metric on the $D^{k+1}$ factor, and a round metric of some small radius $\rho$ on the $S^k$ factor.

Now consider removing a small neighbourhood of $\{0\} \times S^k \subset D^{k+1} \times S^k$ from within the region where the metric is a product. Using a different trivialisation for the surgery then amounts to gluing back this smaller copy of $D^{k+1} \times S^k$ using a boundary diffeomorphism which is different from the identity. The diffeomorphism in question will of course depend on the trivialisation being used.

There is, however, an alternative viewpoint from which we can look at this construction. Consider a neighbourhood $\mathcal N$ of the small copy of $D^{k+1} \times S^k$ which has undergone the twisted gluing into the ambient manifold, that is, consider the small $D^{k+1} \times S^k$ together with a collar. Assume this collar is so small that it lies within the region where the original cap metric is a product. Now imagine cutting out $\mathcal N$, and then gluing it back using the identity map on the boundary. This clearly gives the same topological result as obtained by our previous `twisted' gluing. We can understand this new viewpoint by regarding $\mathcal N$ as a non-trivial $S^k$-bundle over $D^{k+1}$ with a fixed trivialisation near the boundary. To make sense of this, note that $D^{k+1} \times S^k$ is only a trivial bundle over $D^{k+1}$ if we equip it with a globally defined trivialisation. However, if we fix the standard trivialisation on a neighbourhood of the boundary, we can equip $D^{k+1} \times S^k$ with non-trivial bundle structures, that is, such that the boundary trivialisation cannot be extended to a diffeomorphism on all of $D^{k+1} \times S^k.$

We now choose a submersion metric on $\mathcal N$ with base $D^{k+1},$ for which the base and fibre metrics agree with the original choices (i.e. base metric $du^2+\sin(u)ds^2_{k-1}$ and fibre metric $\rho^2ds^2_k$), which glues smoothly into the ambient manifold. Thus near the boundary we need the metric to be a product, however since we are viewing $\mathcal N$ as a non-trivial $S^k$-bundle with fixed boundary trivialisation, it is clear that the submersion metric cannot globally have a product structure. In other words, our choice of normal bundle trivialisation for the surgery results in us having to choose a non-flat connection for our submersion metric, though we are free to choose this so as to be flat in a neighbourhood of the boundary, thus guaranteeing the desired product structure in this region.

It now follows from \cite[Theorem 0.3]{Wr2} that if $\rho$ is chosen sufficiently small, the Ricci curvature of the resulting submersion metric will be strictly positive, irrespective of the connection used. More precisely, there is a constant $\bar{\kappa}>0$ depending on the choice of connection, such that for any $\rho \in (0,\min\{\bar{\rho},\bar{\kappa}\}),$ the submersion metric will be Ricci positive. Thus in the case of a non-standard trivialisation, we must factor the upper bound $\bar{\kappa}$ into our choice of $\rho$ after Corollary \ref{cap_off}.

Exactly the same arguments apply to the surgery on $Z$ in Proposition \ref{V_surgery}, showing that this surgery can also be performed within (local) Ricci positivity for any choice of normal bundle trivialisation, provided our choice of $\rho_0$ is sufficiently small.

Given that we can perform surgery within global $Ric_2>0$ on both $Z$ and $W$ with any normal bundle trivialisations, this completes the proof of Theorem \ref{thm:main}.
\bigskip\bigskip

\appendix

\section{Appendix}

In this appendix we outline the derivation of three sets of Ricci curvature formulas which appear in the main body of the paper without justification. These are the expressions in Lemma \ref{ambient_curvature}, Proposition \ref{ric(g_phi)}, and Lemma \ref{end_ric}. We will take each set of formulas separately, however all rely on the following result. This appeared in \cite{Wr_thesis}, where each expression is carefully derived (using the Kozsul formula).

\begin{lemma}
\label{thesis}
Suppose that $g$ is a submersion metric on the total space of a Riemannian submersion with totally geodesic fibres. Let $\hat{r}$ denote the Ricci curvature of the fibre metric, $\check{r}$ the Ricci curvature of the base metric, let $U,V$ be mutually orthogonal unit vectors with respect to the fibre metric, and $X,Y$ unit vectors with respect to the base metric. Suppose that the fibre metric is then rescaled by a function $\mu^2$ defined on the base. Then if the fibre has dimension $p$, the Ricci curvatures of the resulting metric are as follows.
\begin{align*}
r_{\mu^2}(U/\mu,U/\mu)&=\frac{1}{\mu^2}\hat{r}(U,U)-(p-1)\frac{\|\nabla \mu\|^2}{\mu^2}-\frac{\Delta \mu}{\mu}+\mu^2(AU,AU);\\
r_{\mu^2}(U/\mu,V/\mu)&=\frac{1}{\mu^2}\hat{r}(U,V)+\mu^2(AU,AV);\\
r_{\mu^2}(X,Y)&=\check{r}(X,Y)-2\mu^2(A_X,A_Y)-\frac{p}{\mu}Hess_\mu(X,Y);\\
r_{\mu^2}(X,U/\mu)&=(p+2)(A_XU,\nabla \mu)-\mu(\check{\delta}A(X),U).\\
\end{align*}
Here, $r_{\mu^2}$ denotes the Ricci curvature of the scaled metric, and the $A$-tensor terms are the usual terms for the unscaled metric. All other mixed Ricci curvature terms vanish.
\end{lemma}
\bigskip

\centerline{\bf The formulas of Proposition \ref{ric(g_phi)}.}
\bigskip

The formulas in this Proposition are the Ricci curvatures of the metric $g_\phi$. This was built in two stages: firstly one constructs a submersion metric on a normal annulus bundle to $W$ with base $(W,\check{g})$, fibre $$\Bigl(A^{k+1},d\tau^2+(1-\phi/R)^2dr^2+\bar{\psi}(r,\tau)^2ds^2_{k-1}\Bigr),$$ and connection $\nabla.$ 
Secondly, we rescale the base metric by a function $f^2(\tau)$. For our current purposes however, it will be convenient to view this metric in a different way.

Observe that we can also think of $g_\phi$ as a submersion metric with base $$\Bigl( [0,1] \times (0,R\pi)  \times W\,;\, d\tau^2+\bigl(1-\frac{\phi(\tau)}{R}\bigr)^2dr^2+f^2(\tau)\check{g}\Bigr),$$ fibre $(S^{k-1};\bar{\psi}(r,\tau)^2 ds^2_{k-1}),$ and connection $\bar{\nabla}$, where $\bar{\nabla}$ is the trivial extension of $\nabla$ in $r$ and $\tau$ directions. This change of viewpoint is possible since the $r$ and $\tau$ directions split off: all twisting, both topological and metric is `contained' in the comparison bundle ${\mathcal B}$.

With this new viewpoint, we can use Lemma \ref{thesis} to compute the Ricci curvatures. In order to use the Lemma, we first need to compute the Ricci curvatures of the base, the gradient, Hessian and Laplacian of $\bar{\psi}$ on the base, and the $A$-tensor terms.

We first deal with the base curvatures. 
By computing Christoffel symbols, a straightforward calculation reveals the Ricci curvatures of $$\Bigl({\mathbb R}^2;d\tau^2+(1-\phi(\tau)/R)^2 dr^2\Bigr)$$ to be
\begin{align*}
Ric\Bigl(\frac{1}{1-\phi/R}\frac{\partial}{\partial r},\frac{1}{1-\phi/R}\frac{\partial}{\partial r}\Bigr)&=\frac{\phi''}{R-\phi};\\
Ric\Bigl(\frac{\partial}{\partial \tau},\frac{\partial}{\partial \tau}\Bigr)&=\frac{\phi''}{R-\phi};\\
Ric\Bigl(\frac{\partial}{\partial r},\frac{\partial}{\partial \tau}\Bigr)&=0.\\
\end{align*}
(For convenience we have extended $\phi$ here in the obvious way to a function $\phi:{\mathbb R} \to {\mathbb R}.$ Below, it will also be convenient to view $f$ as a function $f:{\mathbb R} \to {\mathbb R}$, restricting as appropriate on $[0,1].$)

Following on from this, the warped product formulas in \cite[9J]{Be} allow us to compute the Ricci curvature of the base metric $d\tau^2+(1-\phi/R)^2dr^2+f^2\check{g}$ by viewing this as a warped product with fibre $(W,\check{g})$ (which is of course isometric to $(S^k,ds^2_k)$), base $\bigl({\mathbb R}^2;d\tau^2+(1-\phi(\tau)/R)^2 dr^2\bigr),$ and warping function $f^2$. For this, we need the gradient, Hessian and Laplacian of $f$ on the $(r,\tau)$-plane. A simple calculation shows these to be:
\begin{align*}
&\nabla f=f'\frac{\partial}{\partial \tau}; \\
&\text{Hess}_f\bigl(\frac{\partial}{\partial \tau},\frac{\partial}{\partial \tau}\bigr)=f'' ;\\
&\text{Hess}_f\bigl(\frac{\partial}{\partial \tau},\frac{\partial}{\partial r}\bigr)=0; \\
&\text{Hess}_f\Bigl(\frac{1}{1-R^{-1}\phi}\frac{\partial}{\partial r},\frac{1}{1-R^{-1}\phi}\frac{\partial}{\partial r}\Bigr)=-f'\frac{\phi'}{R-\phi}; \\
&\Delta f=f''-\frac{\phi'}{R-\phi}f' .\\
\end{align*}
Using these expressions in the the warped product formulas from \cite[9J]{Be} (bearing in mind that the sign convention for the Laplacian in \cite{Be} is opposite to ours) we immediately obtain the following base-direction Ricci curvatures:
\begin{align*}
&Ric(X_i/f,X_i/f)=\frac{k-1}{f^2}(1-f'^2)-\frac{f''}{f}+\frac{\phi'}{R-\phi}\frac{f'}{f}; \\
&Ric(\partial/\partial \tau,\partial/\partial \tau)=\frac{\phi''}{R-\phi}-k\frac{f''}{f}; \\
&Ric\Bigl(\frac{1}{1-R^{-1}\phi}\frac{\partial}{\partial r},\frac{1}{1-R^{-1}\phi}\frac{\partial}{\partial r}\Bigr)=\frac{\phi''}{R-\phi}+k\frac{\phi'}{R-\phi}\frac{f'}{f}; \\ \\
&Ric(X_i/f,\partial/\partial \tau)=Ric(X_i/f,\partial/\partial r)=Ric(\partial/\partial \tau.\partial/\partial r)=0. \\
\end{align*}

We now turn our attention to the $A$-tensor terms. We start with the $A$-tensor for the comparison bundle $({\mathcal B},g_B).$ We will denote these terms by $A$ as usual. Recall that the base of $\mathcal B$ is $W \cong S^k$, and that the metric $g_B:=\nu(\check g, ds^2_{k-1},\nabla).$ Initially, we consider the metric $d\tau^2+(1-R^{-1}\phi(\tau))^2dr^2+\check{g}$ on ${\mathbb R}^2 \times W$. We then define a new metric on ${\mathbb R}^2 \times {\mathcal B}$ by $$\bar{g}_B:=\nu\Bigl(d\tau^2+(1-R^{-1}\phi)^2 dr^2+\check{g},\ ds^2_{k-1},\ \bar{\nabla}\Bigr).$$ Denote the $A$-tensor of this new metric by $\bar{A}$. As $\bar{g}_B$ is simply the product metric $\bigl(d\tau^2+(1-R^{-1}\phi)^2dr^2\bigr)+g_B$, we see immediately that
$$\bar{A}_{\partial/\partial r} \bullet \equiv \bar{A}_{\partial/\partial \tau} \bullet \equiv 0,$$ and by symmetry
$$\bar{A}_\bullet {\partial/\partial r} \equiv \bar{A}_\bullet {\partial/\partial \tau} \equiv 0,$$ with $\bar{A}=A$ otherwise (i.e. when evaluating on vectors tangent to $\mathcal B$). Similarly $\check{\delta}\bar{A}=\check{\delta}A$ restricted to $T{\mathcal B}.$

Next we introduce the scaling function $f(\tau)$ for the $W$ factor in the base. Thus we consider a new metric $$\bar{g}^f_B:=\nu\Bigl(d\tau^2+(1-R^{-1}\phi)^2 dr^2+f^2(\tau)\check{g},\ ds^2_{k-1},\ \bar{\nabla}\Bigr).$$ A term-by-term computation via the Koszul formula (analogous to that carried out in \cite{Wr_thesis} to prove Lemma \ref{thesis}) yields the following, where $\bar{A}^f$ denotes the $A$-tensor terms for $\bar{g}_B^f.$ 
\begin{align*}
&\bar{A}^f_H H'=\bar{A}_H H' \text{ for any horizontal vectors }H,H'; \\
&\bar{A}^f_{X_i} V=\frac{1}{f^2}\bar{A}_{X_i}V; \\
&\bar{A}^f_{\partial/\partial \tau} \bullet\equiv 0; \\
&\bar{A}^f_{\partial/\partial r} \bullet\equiv 0; \\
&(\check{\delta}\bar{A}^f(X_i),V)_{\bar{g}^f_B}=\frac{1}{f^2}(\check{\delta}\bar{A}(X_i),V)_{\bar{g}_B} \\
&\check{\delta}\bar{A}^f(\partial/\partial \tau)=\check{\delta}\bar{A}^f(\partial/\partial r)=0. \\
\end{align*}

In essence, the above expressions show that the $\bar{A}^f$ terms reduce to scalings of the $A$-tensor terms for $({\mathcal B},g_B)$ for vectors tangent to this bundle, and vanish otherwise. This is consistent with the observation that all the twisting in our extended bundle actually occurs within ${\mathcal B}$.

It follows immediately from the above that the metric $\bar{g}_B^f$ also gives the following terms which appear in the formulas in Lemma \ref{thesis}:
\begin{align*}
&\bigl(\bar{A}^fU,\bar{A}^f V\bigr)_{\bar{g}_B^f}=\frac{1}{f^4}(AU,AV)_{g_B}; \\
&\bigl(\bar{A}^f_X,\bar{A}^f_Y\bigr)_{\bar{g}^f_B}=\frac{1}{f^2}(A_X,A_Y)_{g_B}; \\
&(\check{\delta}\bar{A}^f(X_i),V)_{\bar{g}^f_B}=\frac{1}{f^2}(\check{\delta}A(X_i),V)_{g_B}; \\
& \bar{A}_{X_i} V=\frac{1}{f^2}A_{X_i} V. \\
\end{align*}

Our last task is to introduce the scaling function $\bar{\psi}(r,\tau)$ into the $S^{k-1}$ fibres of our bundle, resulting in the desired metric $g_\phi$. By Lemma \ref{thesis} again, we need to know the gradient, Hessian and Laplacian of $\bar{\psi}$ on the base $$\bigl([0,1]) \times (0,R\pi) \times W;d\tau^2+(1-R^{-1}\phi)^2 dr^2+f^2(\tau)\check{g}\bigr).$$

A computation of Christoffel symbols establishes the following covariant derivatives for the above base metric (with covariant derivatives not listed below equal to 0);
\begin{align*}
&D_{\partial/\partial r} \partial/\partial r=\frac{R-\phi}{R^2}\phi' \frac{\partial}{\partial \tau}; \\
&D_{\partial/\partial r} \partial/\partial \tau=-\frac{\phi'}{R-\phi}\frac{\partial}{\partial r}; \\
&D_{X_i} X_j=-ff'\delta_{ij}\frac{\partial}{\partial \tau}; \\
&D_{X_i} \partial/\partial \tau=D_{\partial/\partial \tau} X_i=\frac{f'}{f}X_i. \\
\end{align*}
The above covariant derivatives then feed in to the computation of the Hessian formulas below.
\begin{align*}
&\text{Hess}_{\bar{\psi}}(\partial/\partial \tau,\partial/\partial \tau)=\bar{\psi}_{\tau \tau}; \\
&\text{Hess}_{\bar{\psi}}(\partial/\partial r,\partial/\partial r)=\bar{\psi}_{rr}-\frac{R-\phi}{R^2}\phi'\bar{\psi}_\tau; \\
&\text{Hess}_{\bar{\psi}}(\partial/\partial \tau,\partial/\partial r)=\bar{\psi}_{\tau r}+\frac{\phi'}{R-\phi}\bar{\psi}_r; \\
&\text{Hess}_{\bar{\psi}}(X_i,X_j)=ff'\delta_{ij}\bar{\psi}_\tau; \\
&\text{Hess}_{\bar{\psi}}(X_i,\partial/\partial \tau)=\text{Hess}_{\bar{\psi}}(X_i,\partial/\partial r)=0. \\
\end{align*}
In particular we have $$\text{Hess}_{\bar{\psi}}\Bigl(\frac{1}{1-R^{-1}\phi}\frac{\partial}{\partial r},\frac{1}{1-R^{-1}\phi}\frac{\partial}{\partial r}\Bigr)=\frac{1}{(1-R^{-1}\phi)^2}\bar{\psi}_{rr}-\frac{\phi'}{R-\phi}\bar{\psi}_\tau;$$
$$\text{Hess}_{\bar{\psi}}\Bigl(\frac{\partial}{\partial \tau},\frac{1}{1-R^{-1}\phi}\frac{\partial}{\partial r}\Bigr)=\frac{1}{1-R^{-1}\phi}\bar{\psi}_{\tau r}+\frac{R\phi'}{(R-\phi)^2}\bar{\psi}_r;$$
$$\text{Hess}_{\bar{\psi}}(X_i/f,X_i/f)=\frac{f'}{f}\bar{\psi}_\tau.$$
We also have $$\|\nabla \bar{\psi}\|^2=\bar{\psi}_\tau^2+\frac{1}{(1-R^{-1}\phi)^2}\bar{\psi}^2_r;$$
$$\Delta \bar{\psi}=\bar{\psi}_{\tau \tau}+\frac{1}{(1-R^{-1}\phi)^2}\bar{\psi}_{rr}+\Bigl(k\frac{f'}{f}-\frac{\phi'}{R-\phi}\Bigr)\bar{\psi}_\tau.$$
Putting the above expressions involving $\bar{\psi}$ together with the previous base curvature and $A$-tensor expressions into the formulas in Lemma \ref{thesis} finally gives the Ricci curvature expressions for $g_\phi$ appearing in Proposition \ref{ric(g_phi)}.
\bigskip

\centerline{\bf The formulas of Lemma \ref{ambient_curvature}.}
\bigskip

This Lemma presents the Ricci curvature formulas for the metric $g$ in terms of the $A$-tensor terms for $({\mathcal B},g_B),$ 
where $g_B=\nu(\check{g},ds^2_{k-1},\nabla).$ 

Consider the manifold $({\mathbb R} \times {\mathcal B}, ds^2+g_B).$ We will view this as a submersion with base ${\mathbb R} \times W$, so $ds^2+g_B=\nu(ds^2+\check{g},ds^2_{k-1},\nabla'),$ where $\nabla'$ is the trivial extension of $\nabla$ in the $\mathbb R$-direction. It is easily checked that the $A$-tensor terms for this submersion are given by $A'_{\partial/\partial s} \bullet \equiv A'_\bullet \partial/\partial s \equiv 0,$ $\check{\delta}A'(\partial/\partial s)=0,$ and $A'=A$ otherwise, where as usual $A$ represents that $A$-tensor terms of $({\mathcal B},g_B).$

Recall that $g=\nu(\check{g},ds^2+\psi^2(s)ds^2_{k-1},\nabla)$ with $s \in [0,\Lambda],$ which can be re-written as $\nu(ds^2+\check{g},\psi^2(s)ds^2_{k-1},\nabla')$, as the $s$-direction splits off and can be `absorbed' into the base. Thus $g$ is obtained by rescaling the fibres of the totally geodesic Riemannian submersion $([0,\Lambda] \times {\mathcal B},ds^2+g_B)$ by $\psi^2(s)$, and of course the Ricci curvatures resulting from this rescaling are described by Lemma \ref{thesis}. An easy computation shows that the gradient, Hessian and Laplacian of $\psi$ on the base are as follows:
\begin{align*}
&\nabla \psi=\psi'\frac{\partial}{\partial s}; \\
&\text{Hess}_\psi\Bigl(\frac{\partial}{\partial s},\frac{\partial}{\partial s}\Bigr)=\psi'',\text{ and Hess}_\psi(\bullet,\bullet)=0 \text{ otherwise}; \\
&\Delta \psi=\psi''. \\
\end{align*}
Using the above data in Lemma \ref{thesis} then produces the desired Ricci curvature formulas for $g$.
\bigskip

\centerline{\bf The formulas of Lemma \ref{end_ric}.}
\bigskip

This Lemma presents the Ricci curvature expressions for the metric $$du^2+\nu\bigl(F^2(u)ds^2_k,h^2(u)\omega, \nabla\bigr).$$ 
It will be convenient to view this metric as 
\begin{equation}\label{no1}
\nu(du^2+F^2(u)ds^2_k,h^2(u)\omega,\nabla'),
\end{equation}
where $\nabla'$ is the trivial extension of $\nabla$ in the $u$-direction. In order to apply Lemma \ref{thesis} to this metric, we first need to consider the metric
\begin{equation}\label{no2}
\nu(du^2+F^2(u)ds^2_k,\omega,\nabla').
\end{equation}
Once we have established the curvatures of (\ref{no2}), Lemma \ref{thesis} then shows how to incorporate the fibre scalings by $h^2(u).$

To compute the Ricci curvatures of (\ref{no2}), we need the Ricci curvatures of the base and fibres, and the $A$-tensor terms. We also need the gradient, Hessian and Laplacian of $h$ on the base $({\mathbb R} \times S^k,du^2+F^2(u)ds^2_k).$ As the metrics on both base and fibre are (single) warped products, we can simply write down these expressions in terms of $F$, $\bar{\theta}$, and their derivatives. Similarly, the gradient, Hessian and Laplacian of $h$ are easily written down, and yield analogous expressions to computations elsewhere in this paper. The non-trivial part is understanding the $A$-tensor.

Our approach to identifying the $A$-tensor for (\ref{no2}) is to first observe that this metric can be re-written in the form
\begin{equation}\label{no3}
\nu(d\eta^2+du^2+F^2(u)ds^2_k,\bar{\theta}^2(\eta)ds^2_{k-1},\bar{\nabla}),
\end{equation}
where $\bar{\nabla}$ is the trivial extension of $\nabla$ in both $u$ and $\eta$-directions. The reason this alternative description exists is that without the $h$-scaling on the fibres, the $\eta$-direction in effect splits off and can be absorbed into the base.

Ultimately we want to express the $A$-tensor for (\ref{no2}), $A'$, in terms of the $A$-tensor for $({\mathcal B},g_B),$ which as usual we will simply denote by $A$. However we begin by comparing $A'$ with the $A$-tensor for (\ref{no3}), which we will label $\bar{A}$.

As the metrics (\ref{no2}) and (\ref{no3}) are identical, even though the submersion structures and $A$-tensors differ, the covariant derivatives $D$ are the same. Therefore to compare $A'$ and $\bar{A}$, we only have to keep track of vectors in the $\eta$-direction, 
since $A'$ and $\bar{A}$ only differ in their interpretation of horizontal and vertical vectors, and the $\eta$-direction is precisely the difference.

A straightforward computation, for example using the Koszul formula, shows that the $A'$-tensor terms we need for the curvature formulas are either the same as the corresponding $\bar{A}$-expression when $\partial/\partial\eta$ is not present, or equal to zero when $A'$-term is being evaluated on $\partial/\partial\eta.$

The next step is to use $\bar{A}$ to calculate the $A$-tensor, $\tilde{A}$, of the metric
\begin{equation}\label{no4}
\nu(d\eta^2+du^2+F^2(u)ds^2_k,ds^2_{k-1},\bar{\nabla}).
\end{equation}
Here we essentially want to use the formulas of Lemma \ref{thesis} in reverse. In practice, this involves an elementary calculation using the $A$-tensor transformation formulas established in \cite{Wr_thesis} (which are used to prove Lemma \ref{thesis}). This yields, for example, $(\bar{A}V,\bar{A}V)=\bar{\theta}^4(\tilde{A}V,\tilde{A}V),$ $(\bar{A}_X,\bar{A}_Y)=\bar{\theta}^2(\tilde{A}_X,\tilde{A}_Y),$ and so on. (We have suppressed the subscripts here from the bracketed expressions which indicate the specific metric to be used.)

Finally, we can express the $\tilde{A}$-terms in terms of the $A$-tensor terms for $({\mathcal B},g_B)$ using the transformation formulas displayed earlier in this Appendix in relation to Proposition \ref{ric(g_phi)}, where we replace $\tau$ by $u$, $r$ by $\eta$, and set $\phi \equiv 0.$

\bigskip\bigskip

\noindent{\it Diarmuid Crowley, 
School of Mathematics \& Statistics,
The University of Melbourne,
Parkville, VIC, 3010,
Australia.
Email: dcrowley@unimelb.edu.au}\\

\medskip
\noindent {\it David Wraith, 
Department of Mathematics and Statistics, 
National University of Ireland Maynooth, Maynooth, 
County Kildare, 
Ireland. 
Email: david.wraith@mu.ie.}

\end{document}